\documentclass{gtart}

\def\ifplaintex{\expandafter\ifx\csname documentclass\endcsname\relax}

\ifplaintex 
\else
\fi

\expandafter\ifx\csname beginpicture\endcsname\relax
\expandafter\ifx\csname documentclass\endcsname\relax
\input pictex \else
\input prepictex \input pictex \input postpictex \fi\fi

\def\gt{{\mathsurround=0pt\it $\cal G\mskip-2mu$eometry \&\ 
$\cal T\!\!$opology}}        %  journal title in recommended style

\def\gtp{{\mathsurround=0pt\it $\cal G\mskip-2mu$eometry \&\ 
$\cal T\!\!$opology $\cal P\!$ublications}}  % GT publications

\def\lognumber#1{\def\thelognumber{#1}}
\def\volumenumber#1{\def\thevolumenumber{#1}}
\def\papernumber#1{\def\thepapernumber{#1}}
\def\volumeyear#1{\def\thevolumeyear{#1}}

\def\pagenumbers#1#2{\def\startpage{#1}\def\finishpage{#2}}
\def\published#1{\def\publishdate{#1}}
\def\proposed#1{\def\theproposer{#1}}
\def\seconded#1{\def\theseconders{#1}}
\def\received#1{\def\receiveddate{#1}}

\def\accepted#1{\def\accepteddate{#1}}
\def\asciititle#1{\def\theasciititle{#1}}

\def\asciiaddress#1{\def\theasciiaddress{#1}}
\def\asciiemail#1{\def\theasciiemail{#1}}
\long\def\asciiabstract#1{\long\def\theasciiabstract{#1}}
\def\asciikeywords#1{\def\theasciikeywords{#1}}
\def\shortauthors#1{\def\theshortauthors{#1}}

\let\\\par\let\thelognumber\relax
\let\thevolumenumber\relax\let\thepapernumber\relax
\let\thevolumeyear\relax\let\thesamplenumber\relax\let\startpage\relax
\let\finishpage\relax\let\publishdate\relax\let\receiveddate\relax
\let\reviseddate\relax\let\accepteddate\relax\let\theasciititle\relax
\let\theasciiauthors\relax\let\theasciiaddress\relax
\let\theasciiabstract\relax\let\theasciikeywords\relax
\let\theasciiemail\relax\let\theshortauthors\relax\let\theshorttitle\relax

\long\def\maketitlep{   % start of definition of \maketitlep

\count0=\startpage

\gt\hfill      %   Journal title (top left) 
\beginpicture
\setcoordinatesystem units <0.33truein, 0.33truein> point at 2.2 0.9
\setplotsymbol ({$\cal G$})
\plotsymbolspacing=9truept
\circulararc 315 degrees from 0 1 center at 0 0
\setplotsymbol ({$\cal T$})
\circulararc 315 degrees from 1 -1 center at 1 0
\endpicture
\break
{\small\ifx\thesamplenumber\relax % sample?  
Volume \else Sample
\fi\thevolumenumber\ (\thevolumeyear)
\startpage--\finishpage\nl
Published: \publishdate}
\vglue 0.5truein plus 0.4fil minus 0.1truein

{\parskip=0pt\leftskip 0pt plus 1fil\def\\{\par\smallskip}{\ifplaintex\large
\else\Large\fi\bf\thetitle}\par\medskip}   

\vglue 0pt plus 0.1fil 

{\parskip=0pt\leftskip 0pt plus 1fil\def\\{\par}{\sc\theauthors}
\par\medskip}

\vglue 0pt plus 0.1fil 

{\small\parskip=0pt\let\newline\\
{\leftskip 0pt plus 1fil\def\\{\par}{\sl\theaddress}\par}
\expandafter\ifx\theemail\relax    % email address?
\relax\else\vglue 5pt plus 0.02fil minus 2pt\def\\{\stdspace{\rm 
and}\stdspace} 
\cl{Email:\stdspace\tt\theemail}\fi
\ifx\theurl\relax                  % URL given?
\relax\else\vglue 5pt plus 0.02fil minus 2pt\def\\{\stdspace{\rm 
and}\stdspace}
\cl{URL:\stdspace\tt\theurl}\fi\par}

\vglue 7pt plus 0.3fil minus 3pt

{\bf Abstract}
\vglue 5pt plus 0.1fil minus 2pt

\theabstract

\vglue 7pt plus 0.3fil minus 3pt

{\bf AMS Classification numbers}\quad Primary:\quad \theprimaryclass

Secondary:\quad \thesecondaryclass

\vglue 5pt plus 0.3fil minus 2pt

{\bf Keywords}\quad \thekeywords

\vglue 10pt plus 0.5fil minus 5pt

{\small  Proposed: \theproposer\hfill Received: \receiveddate\nl
Seconded: \theseconders\hfill 
\ifx\reviseddate\relax                         % paper revised?
Accepted: \accepteddate                        % no
\else
Revised: \reviseddate                          % yes
\fi}
\eject
}       %  end of definition of \maketitlep

\let\maketitlepage\maketitlep
\let\maketitle\maketitlepage

\font\phead=cmsl9 scaled 950
\font=cmsl9 scaled 1050
\font\pnum=cmbx10 scaled 913
\font\lnum=cmbx10 
\font\pfoot=cmsl9 scaled 950
\font=cmsl9 scaled 1050
\ifplaintex
\headline{\vbox to 0pt{\vskip -4.5mm\line{\small\phead\ifnum
\count0=\startpage ISSN 1364-0380 (on line)
1465-3060 (printed) \hfill {\pnum\folio}\else\ifodd\count0\def\\{ }% 
\ifx\theshorttitle\relax\thetitle\else\theshorttitle\fi\hfill{\pnum\folio}
\else\def\\{ and }{\pnum\folio}\hfill\ifx\theshortauthors\relax\theauthors
\else\theshortauthors\fi\fi\fi}\vss}}
\footline{\vbox to 0pt{\vglue 0mm\line{\small\pfoot\ifnum\count0=\startpage
\copyright\ \gtp\hfill\else
\gt, Volume \thevolumenumber\ (\thevolumeyear)\hfill\fi}\vss
}}
\else
\makeatletter
\def\@oddhead{{\small\lhead\ifnum\count0=\startpage ISSN 1364-0380 (on line)
1465-3060 (printed) \hfill {\lnum\number\count0}\else\ifodd\count0
\def\\{ }\ifx\theshorttitle\relax \thetitle \else\theshorttitle\fi\hfill
{\lnum\number\count0}\else\def\\{ and }{\lnum\number\count0}
\hfill\ifx\theshortauthors\relax 
\theauthors\else\theshortauthors\fi\fi\fi}}\def\@evenhead{\@oddhead}
\def\@oddfoot{\small\lfoot\ifnum\count0=\startpage\copyright\ \gtp\hfill\else
\gt, Volume \thevolumenumber\ (\thevolumeyear)\hfill\fi}
\def\@evenfoot{\@oddfoot}
\makeatother
\fi

\newwrite\gtoutfile
\long\gdef\makeheadfile{  %%% start of definition of \makeheadfile
{\def\\{, }\def\s{ }
\immediate\openout\gtoutfile head.xxx
\immediate\write\gtoutfile{To: math@arxiv.org}
\immediate\write\gtoutfile{Subject: put or rep NNNNN:pppp}
\immediate\write\gtoutfile{--text follows this line--}
\immediate\write\gtoutfile{Proxy-for: \ifx\theasciiauthors\relax
\theauthors\else\theasciiauthors\fi\s<\ifx\theasciiemail\relax\theemail\else\theasciiemail\fi>}
\immediate\write\gtoutfile{\noexpand\\}
\immediate\write\gtoutfile{Authors: \ifx\theasciiauthors\relax
\theauthors\else\theasciiauthors\fi}
\immediate\write\gtoutfile{Title: \ifx\theasciititle\relax
\thetitle\else\theasciititle\fi}
\immediate\write\gtoutfile{Subj-class: GT or SG or MG etc}
\immediate\write\gtoutfile{MSC-class: \theprimaryclass\ifx\thesecondaryclass\relax\else, \thesecondaryclass\fi}
\immediate\write\gtoutfile{Journal-ref: Geom. Topol. \thevolumenumber
(\thevolumeyear) \startpage-\finishpage}
\immediate\write\gtoutfile{Comments: Published by Geometry and Topology at}
\immediate\write\gtoutfile{\s\s http://www.maths.warwick.ac.uk/gt/GTVol\thevolumenumber/paper\thepapernumber.abs.html}
\immediate\write\gtoutfile{\noexpand\\}
\immediate\write\gtoutfile{}
\ifx\theasciiabstract\relax
\immediate\write\gtoutfile{\theabstract}\else
\immediate\write\gtoutfile{\theasciiabstract}\fi
\immediate\write\gtoutfile{}
\immediate\write\gtoutfile{\noexpand\\}
\immediate\write\gtoutfile{}
\immediate\closeout\gtoutfile}}  %%% end of definition of \makeheadfile

\def\maketitlepage{\maketitlep\makeheadfile}
\let\maketitle\maketitlepage

\def\ifplaintex{\expandafter\ifx\csname documentclass\endcsname\relax}

\ifplaintex 
\else
\fi

\expandafter\ifx\csname beginpicture\endcsname\relax
\expandafter\ifx\csname documentclass\endcsname\relax
\input pictex \else
\input prepictex \input pictex \input postpictex \fi\fi

\def\gt{{\mathsurround=0pt\it $\cal G\mskip-2mu$eometry \&\ 
$\cal T\!\!$opology}}        %  journal title in recommended style

\def\gtp{{\mathsurround=0pt\it $\cal G\mskip-2mu$eometry \&\ 
$\cal T\!\!$opology $\cal P\!$ublications}}  % GT publications

\def\lognumber#1{\def\thelognumber{#1}}
\def\volumenumber#1{\def\thevolumenumber{#1}}
\def\papernumber#1{\def\thepapernumber{#1}}
\def\volumeyear#1{\def\thevolumeyear{#1}}

\def\pagenumbers#1#2{\def\startpage{#1}\def\finishpage{#2}}
\def\published#1{\def\publishdate{#1}}
\def\proposed#1{\def\theproposer{#1}}
\def\seconded#1{\def\theseconders{#1}}
\def\received#1{\def\receiveddate{#1}}

\def\accepted#1{\def\accepteddate{#1}}
\def\asciititle#1{\def\theasciititle{#1}}

\def\asciiaddress#1{\def\theasciiaddress{#1}}
\def\asciiemail#1{\def\theasciiemail{#1}}
\long\def\asciiabstract#1{\long\def\theasciiabstract{#1}}
\def\asciikeywords#1{\def\theasciikeywords{#1}}
\def\shortauthors#1{\def\theshortauthors{#1}}

\let\\\par\let\thelognumber\relax
\let\thevolumenumber\relax\let\thepapernumber\relax
\let\thevolumeyear\relax\let\thesamplenumber\relax\let\startpage\relax
\let\finishpage\relax\let\publishdate\relax\let\receiveddate\relax
\let\reviseddate\relax\let\accepteddate\relax\let\theasciititle\relax
\let\theasciiauthors\relax\let\theasciiaddress\relax
\let\theasciiabstract\relax\let\theasciikeywords\relax
\let\theasciiemail\relax\let\theshortauthors\relax\let\theshorttitle\relax

\long\def\maketitlep{   % start of definition of \maketitlep

\count0=\startpage

\gt\hfill      %   Journal title (top left) 
\beginpicture
\setcoordinatesystem units <0.33truein, 0.33truein> point at 2.2 0.9
\setplotsymbol ({$\cal G$})
\plotsymbolspacing=9truept
\circulararc 315 degrees from 0 1 center at 0 0
\setplotsymbol ({$\cal T$})
\circulararc 315 degrees from 1 -1 center at 1 0
\endpicture
\break
{\small\ifx\thesamplenumber\relax % sample?  
Volume \else Sample
\fi\thevolumenumber\ (\thevolumeyear)
\startpage--\finishpage\nl
Published: \publishdate}
\vglue 0.5truein plus 0.4fil minus 0.1truein

{\parskip=0pt\leftskip 0pt plus 1fil\def\\{\par\smallskip}{\ifplaintex\large
\else\Large\fi\bf\thetitle}\par\medskip}   

\vglue 0pt plus 0.1fil 

{\parskip=0pt\leftskip 0pt plus 1fil\def\\{\par}{\sc\theauthors}
\par\medskip}

\vglue 0pt plus 0.1fil 

{\small\parskip=0pt\let\newline\\
{\leftskip 0pt plus 1fil\def\\{\par}{\sl\theaddress}\par}
\expandafter\ifx\theemail\relax    % email address?
\relax\else\vglue 5pt plus 0.02fil minus 2pt\def\\{\stdspace{\rm 
and}\stdspace} 
\cl{Email:\stdspace\tt\theemail}\fi
\ifx\theurl\relax                  % URL given?
\relax\else\vglue 5pt plus 0.02fil minus 2pt\def\\{\stdspace{\rm 
and}\stdspace}
\cl{URL:\stdspace\tt\theurl}\fi\par}

\vglue 7pt plus 0.3fil minus 3pt

{\bf Abstract}
\vglue 5pt plus 0.1fil minus 2pt

\theabstract

\vglue 7pt plus 0.3fil minus 3pt

{\bf AMS Classification numbers}\quad Primary:\quad \theprimaryclass

Secondary:\quad \thesecondaryclass

\vglue 5pt plus 0.3fil minus 2pt

{\bf Keywords}\quad \thekeywords

\vglue 10pt plus 0.5fil minus 5pt

{\small  Proposed: \theproposer\hfill Received: \receiveddate\nl
Seconded: \theseconders\hfill 
\ifx\reviseddate\relax                         % paper revised?
Accepted: \accepteddate                        % no
\else
Revised: \reviseddate                          % yes
\fi}
\eject
}       %  end of definition of \maketitlep

\let\maketitlepage\maketitlep
\let\maketitle\maketitlepage

\font\phead=cmsl9 scaled 950
\font=cmsl9 scaled 1050
\font\pnum=cmbx10 scaled 913
\font\lnum=cmbx10 
\font\pfoot=cmsl9 scaled 950
\font=cmsl9 scaled 1050
\ifplaintex
\headline{\vbox to 0pt{\vskip -4.5mm\line{\small\phead\ifnum
\count0=\startpage ISSN 1364-0380 (on line)
1465-3060 (printed) \hfill {\pnum\folio}\else\ifodd\count0\def\\{ }% 
\ifx\theshorttitle\relax\thetitle\else\theshorttitle\fi\hfill{\pnum\folio}
\else\def\\{ and }{\pnum\folio}\hfill\ifx\theshortauthors\relax\theauthors
\else\theshortauthors\fi\fi\fi}\vss}}
\footline{\vbox to 0pt{\vglue 0mm\line{\small\pfoot\ifnum\count0=\startpage
\copyright\ \gtp\hfill\else
\gt, Volume \thevolumenumber\ (\thevolumeyear)\hfill\fi}\vss
}}
\else
\makeatletter
\def\@oddhead{{\small\lhead\ifnum\count0=\startpage ISSN 1364-0380 (on line)
1465-3060 (printed) \hfill {\lnum\number\count0}\else\ifodd\count0
\def\\{ }\ifx\theshorttitle\relax \thetitle \else\theshorttitle\fi\hfill
{\lnum\number\count0}\else\def\\{ and }{\lnum\number\count0}
\hfill\ifx\theshortauthors\relax 
\theauthors\else\theshortauthors\fi\fi\fi}}\def\@evenhead{\@oddhead}
\def\@oddfoot{\small\lfoot\ifnum\count0=\startpage\copyright\ \gtp\hfill\else
\gt, Volume \thevolumenumber\ (\thevolumeyear)\hfill\fi}
\def\@evenfoot{\@oddfoot}
\makeatother
\fi

\newwrite\gtoutfile
\long\gdef\makeheadfile{  %%% start of definition of \makeheadfile
{\def\\{, }\def\s{ }
\immediate\openout\gtoutfile head.xxx
\immediate\write\gtoutfile{To: math@arxiv.org}
\immediate\write\gtoutfile{Subject: put or rep NNNNN:pppp}
\immediate\write\gtoutfile{--text follows this line--}
\immediate\write\gtoutfile{Proxy-for: \ifx\theasciiauthors\relax
\theauthors\else\theasciiauthors\fi\s<\ifx\theasciiemail\relax\theemail\else\theasciiemail\fi>}
\immediate\write\gtoutfile{\noexpand\\}
\immediate\write\gtoutfile{Authors: \ifx\theasciiauthors\relax
\theauthors\else\theasciiauthors\fi}
\immediate\write\gtoutfile{Title: \ifx\theasciititle\relax
\thetitle\else\theasciititle\fi}
\immediate\write\gtoutfile{Subj-class: GT or SG or MG etc}
\immediate\write\gtoutfile{MSC-class: \theprimaryclass\ifx\thesecondaryclass\relax\else, \thesecondaryclass\fi}
\immediate\write\gtoutfile{Journal-ref: Geom. Topol. \thevolumenumber
(\thevolumeyear) \startpage-\finishpage}
\immediate\write\gtoutfile{Comments: Published by Geometry and Topology at}
\immediate\write\gtoutfile{\s\s http://www.maths.warwick.ac.uk/gt/GTVol\thevolumenumber/paper\thepapernumber.abs.html}
\immediate\write\gtoutfile{\noexpand\\}
\immediate\write\gtoutfile{}
\ifx\theasciiabstract\relax
\immediate\write\gtoutfile{\theabstract}\else
\immediate\write\gtoutfile{\theasciiabstract}\fi
\immediate\write\gtoutfile{}
\immediate\write\gtoutfile{\noexpand\\}
\immediate\write\gtoutfile{}
\immediate\closeout\gtoutfile}}  %%% end of definition of \makeheadfile

\def\maketitlepage{\maketitlep\makeheadfile}
\let\maketitle\maketitlepage

\lognumber{144}
\volumenumber{5}\papernumber{3}\volumeyear{2001}
\pagenumbers{75}{108}
\accepted{28 January 2001}
\proposed{Robion Kirby}
\seconded{Joan Birman, Cameron Gordon}
\received{19 October 2000}
\published{10 February 2001}

\usepackage{amsmath,psfig,amssymb}

\newcommand{\N}{\mathbb N}
\newcommand{\R}{\mathbb R}
\newcommand{\Z}{\mathbb Z}
\newcommand{\oZ}{\otimes\mathbb Z[1/2]}
\newcommand{\GG}{\Gamma}
\newcommand{\GD}{\Delta}
\newcommand{\Gd}{\delta}
\newcommand{\Gs}{\sigma}
\newcommand{\GS}{\Sigma}
\newcommand{\cA}{\mathcal A}
\newcommand{\tA}{\widetilde{\mathcal A}}
\newcommand{\cF}{\mathcal F}
\newcommand{\cFY}{\mathcal F^Y}
\newcommand{\cFa}{\mathcal F^{as}}
\newcommand{\cFb}{\mathcal F^b}
\newcommand{\cFbl}{\mathcal F^{bl}}
\newcommand{\cG}{\mathcal G}

\newcommand{\tpsi}{\widetilde\psi}

\newcommand{\lk}{{\mathrm lk}}
\newcommand{\sminus}{\smallsetminus}

\newtheorem{thm}{Theorem}[section]
\newtheorem{lem}[thm]{Lemma}
\newtheorem{prop}[thm]{Proposition}
\newtheorem{cor}[thm]{Corollary}
\theoremstyle{remark}
\newtheorem{rem}[thm]{Remark}

\def\eq#1{\underset{#1}{=}}
\def\fig#1{\vcenter{\psfig{figure=#1,silent=}}}
\def\centerfig#1{\centerline{\psfig{figure=#1,silent=}}}
\def\fti{finite type invariant}
\def\BQ{\mathbb Q}

\begin{document}
\title{Calculus of clovers and finite type invariants\\of 3--manifolds}
\asciititle{Calculus of clovers and finite type invariants
of 3-manifolds}
\authors{Stavros Garoufalidis\\Mikhail Goussarov\\Michael Polyak}
\shortauthors{Garoufalidis, Goussarov and Polyak}
\address{{\rm SG:}\ \ School of Mathematics, Georgia Institute of
Technology\\Atlanta, GA 30332-0160 USA\\\smallskip\\
{\rm MP:}\ \ School of Mathematics, Tel-Aviv 
University\\69978 Tel-Aviv, Israel}
\asciiaddress{SG: School of Mathematics, Georgia Institute of
Technology\\Atlanta, GA 30332-0160 USA\\MP: School of Mathematics, Tel-Aviv 
University\\69978 Tel-Aviv, Israel}
\email{stavros@math.gatech.edu}
\secondemail{polyak@math.tau.ac.il}
\asciiemail{stavros@math.gatech.edu, polyak@math.tau.ac.il}

\begin{abstract}
A clover is a framed trivalent graph with some additional structure,
embedded in a 3--manifold. We define surgery on clovers, generalizing
surgery on Y--graphs used earlier by the second author to define a new
theory of finite-type invariants of 3--manifolds.
We give a systematic exposition of a topological calculus of clovers
and use it to deduce some important results about the corresponding
theory of finite type invariants.
In particular, we give a description of the weight systems in terms of
uni-trivalent graphs modulo the AS and IHX relations, reminiscent of
the similar results for links. We then compare several definitions of
finite type invariants of homology spheres (based on surgery on
Y--graphs, blinks, algebraically split links, and boundary links) and
prove in a self-contained way their equivalence.
\end{abstract}
\asciiabstract{
A clover is a framed trivalent graph with some additional structure,
embedded in a 3-manifold. We define surgery on clovers, generalizing
surgery on Y-graphs used earlier by the second author to define a new
theory of finite-type invariants of 3--manifolds.
We give a systematic exposition of a topological calculus of clovers
and use it to deduce some important results about the corresponding
theory of finite type invariants.
In particular, we give a description of the weight systems in terms of
uni-trivalent graphs modulo the AS and IHX relations, reminiscent of
the similar results for links. We then compare several definitions of
finite type invariants of homology spheres (based on surgery on
Y--graphs, blinks, algebraically split links, and boundary links) and
prove in a self-contained way their equivalence.}

\asciikeywords{3-manifolds, Y-graphs, finite type invariants, clovers}

\maketitlepage

\section{Introduction}\label{s1}
\subsection{Finite type invariants}
Polynomials play a fundamental role in mathematics.
Here is one of the possible definitions of a polynomial
using a discrete version of the $n$-th derivative.
Let $V$ be a vector space over a field $k$.
For $(x_1,\dots,x_n)\in V^n$ and $\Gs\in\{0,1\}^n$,
denote $x_\Gs=\sum_{i:\Gs_i=1}x_i$ and $|\Gs|=\sum_i\Gs_i$.
A function $f\co V\to k$ is a polynomial of degree less
than $n$, iff $\sum_\Gs(-1)^{|\Gs|}f(x_\Gs)=0$ for any
$(x_1,\dots,x_n)\in V^n$.
This definition has a significant advantage over the
standard definition of a polynomial, since it can be
easily considered in a more general setup.
In particular, $V$ and $k$ may be just abelian groups.

Finite type invariants play the role of polynomial
functions in topology. In this setting, however, the
situation is somewhat more tricky, and involves a
choice of an appropriate analog of both the zero
element $0$ and the addition operation $+$ in $V$.
In the last decade, this general idea was successfully
used in knot theory under the name of Vassiliev
invariants.
This paper is devoted exclusively to a study of a theory
of finite type invariants of 3--manifolds, initiated
by the second author.

The present paper was started in September 1998
in an attempt to present unpublished results by the
second author.
In June 1999, when this paper was still in its initial
stage, the second author passed away in a tragic
accident. As a result, many details and proofs were
lost. We tried our best to reconstruct them and to
present our results in his style and spirit.

The  authors were partially supported by BSF
grant 97-00398 and by NSF grant DMS-98-00703.

\subsection{Finite type invariants of 3--manifolds}
Several authors developed different theories
of finite type invariants of integer homology spheres.
All these approaches are based on the cut-and-paste
technique along different types of handlebodies.
Historically the first, the theory of T Ohtsuki
\cite{Oh} is based on surgery on algebraically split
links and is by now relatively well understood.
The first author \cite{Ga} proposed then a theory
based on surgery on boundary links, which remained
relatively undeveloped.
Subsequently, the first author and J Levine \cite{GL}
considered theories based on surgery on blinks, and on
gluings along surfaces using elements of the Torelli
group.
All these theories turn out to be equivalent
(up to a 2--torsion), see Section \ref{sub_main}.

The second author initiated in \cite{Gu2} a new
theory of finite type invariants for arbitrary
3--manifolds, based on surgery along Y--graphs
(see \cite{Gu1, Gu2}).
The latter theory has some important advantages.
Firstly, it allows a unified treatment of links,
graphs and 3--manifolds (possibly with boundary or
a Spin--structure). It is by now relatively well
understood for rational homology spheres. Secondly,
other theories mentioned above have a technical
drawback: the corresponding classes of surgery
links are not preserved under handle slides.
Finally, this theory comes equipped with a powerful
topological calculus (introduced in \cite{Gu1}),
that is well-suited for a study of 3--manifolds
and for explicit computations.
This calculus is similar to Kirby's calculus of
framed links, but instead of cutting, twisting
and regluing solid tori, in the clover calculus
one performs similar operations with solid
handlebodies of higher genus. The role of framed
links is played here by {\em clovers}, which are
framed trivalent graphs with some additional
structure.

The main goals of the paper are:
\begin{itemize}
\item
To give an exposition of the calculus of clovers
(see Sections \ref{sec_clovers}--\ref{sec_calculus}),
 and to describe the graded spaces of the theory
in terms of uni-trivalent graphs (see Sections
\ref{sub_main} and \ref{sub_graphs}).
\item
To show, as an application,
 in an elementary and self-contained way,
that all the above mentioned theories of \fti s of
integer homology spheres, properly indexed, coincide
over $\mathbb Z[1/2]$ (see Section \ref{sub_main}).
\end{itemize}

Independently and being uninformed about the results
of the second author, K Habiro developed a theory
of claspers \cite{Ha} and used it to study
\fti s of links. He also announced in
\cite[Section 8.4]{Ha} that his theory can be
extended to the study of \fti s of 3--manifolds
 and gave a brief outline of this extension.
It is based on {\em allowable graph claspers,}
which differ from the type of claspers considered
in the rest of Habiro's paper.
In particular, the basic claspers, which were
Habiro's main tool in the theory for links, are
not allowed in this setting.
In general it is plausible that invariants  of
links can be extended to invariants of 3--manifolds.
In the case in hand, however, this extension is
rather nontrivial and involves new ideas presented here.

The initial data, the way to describe the basic
objects, and the graphical calculus in the clasper
and clover theories are somewhat different.
However, while different in some aspects, these
two theories essentially coincide, with surgery
on allowable graph claspers being equivalent to
surgery on clovers.
This independent discovery of almost the same
theory by two non-interacting researchers looks
rather promising, ruling out an element of
arbitrariness.
Some results announced by K Habiro \cite[Section
8.4]{Ha} seem to overlap with those discussed in
Sections \ref{sec_clovers}--\ref{sec_graded} below.
Unfortunately, it seems that presently there is no
available manuscript with the exact statements and
proofs of Habiro's results for 3--manifolds.
We apologize in advance for any possible overlaps.
It should be also mentioned that T Cochran and
P Melvin \cite{CM} proposed an extension of Ohtsuki's
theory to arbitrary 3--manifolds.
Their theory, unlike the one considered in this paper,
preserves triple cup-products, and the relation between
these theories is more complicated.

\subsection{Y--graphs}
Throughout this paper by a manifold we mean a smooth
oriented compact connected 3--manifold.

\begin{figure}[htb]
\centerfig{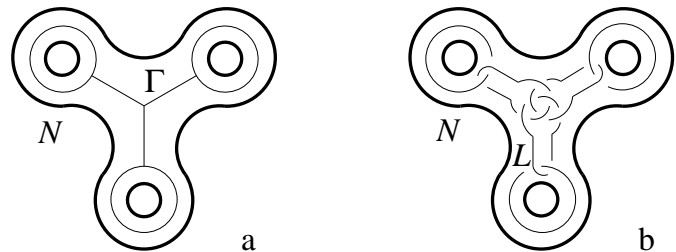,height=1.3in}
\caption{Y--graph and the corresponding surgery link}
\label{fig1}
\end{figure}

We recall some definitions from \cite{Gu2}.
The graph $\GG\subset\R^3$ shown in  Figure \ref{fig1}a
is called a {\em standard Y--graph.}
The edges of $\GG$ are framed with a vector field normal
to the plane of the picture.
A framed graph $G$ in a 3--manifold $M$ is called a {\em
Y--graph\/}, if it is the image of $\GG$ under a smooth
embedding $\phi_G\co N\to M$ of a neighborhood $N$ of $\GG$.
The embedding $\phi_G$ can be recovered from $G$ up to
isotopy.
Let $L$ be the six component link in $N$ shown in  Figure
\ref{fig1}b.
All components of $L$ are 0--framed.
Surgery on $M$ along the framed link $\phi_G(L)$ is
called a {\em surgery along $G$,} or a {\em Y--surgery.}
A Y--surgery can be realized by cutting a solid genus 3
handlebody $\phi_G(N)$ and regluing it in another way
(see \cite{Ma}).
The resulting manifold is defined up to a diffeomorphism,
which is the identity outside of a small neighborhood of
the Y--graph.
Denote it by $M_G$.

An equivalent surgery, under the name of {\em Borromean
surgery} was considered by S~Matveev.
As shown in \cite{Ma}, two manifolds can be connected
by a sequence of Borromean surgeries, if and only if
they have the same homology and the linking pairing in
$\operatorname{Tors}H_1$.
In particular, $M$ can be obtained from $S^3$ by
Y--surgeries, iff $M$ is an integer homology sphere
($\Z$HS in short).

A {\em Y--link} $G$ in a manifold $M$ is a collection
of disjoint Y--graphs in $M$. By surgery on $G$ we
mean surgery on each Y--graph in $G$.
Pushing out the Y--graphs of the following surgeries
from the previously reglued handlebodies, we can present
a sequence of Y--surgeries by surgery on a Y--link.

\subsection{Colored Y--links and the filtration $\cFY_n(M)$}
An {\em $n$--coloring} of a Y--link $G$ is its splitting
into $n$ disjoint subcollections $G_1$, $G_2$,\dots,
$G_n$ of Y--graphs.
The coloring is {\em simple}, if each subcollection
$G_i$ consists of a single Y--graph.

For an $n$--tuple $\Gs=\{\Gs_1,\dots,\Gs_n\}\in\{0,1\}^n$,
let $|\Gs|=\sum_i\Gs_i$ be the number of ones in $\Gs$
and put $G(\Gs)=\bigcup_{i:\Gs_i=1}G_i$.
In an abelian group $\cF$ generated by all 3--manifolds,
put $$[M,G]=\sum_\Gs (-1)^{|\Gs|}M_{G(\Gs)}$$ for a
colored Y--link $G$ in $M$.

Let $\cFY_n(M)$ be the subgroup of $\cF$, generated
by $[M,G]$ for all Y--links $G$, colored in at least
$n$ colors.
This defines a decreasing filtration $\cFY(M)=\cFY_0(M)
\supset\cFY_1(M)\supset\cFY_2(M)\supset\dots$.

Given a Y--link $G$ in $M$ without a specified
coloring, by $[M,G]$ we mean the above alternating
sum for a simple coloring of $G$.

\begin{lem} The subgroup $\cFY_n(M)$ is generated
by $[M,G]$ for all simply colored Y--links $G$ with
at least $n$ components.
\end{lem}
\begin{proof}
Let $G$ be an $n$--colored Y--link with $m$ components.
We call $d=n-m$ the defect of the coloring.
If $d=0$, the coloring is simple.
Suppose $d>0$; then there are at least two components
of the same color $j$.
Therefore, the $j$--colored Y--sublink of $G$ can be
split into two disjoint non-empty parts $G'$ and $G''$.
It is easy to see that
$$[M,G]=[M,G\sminus G']+[M,G\sminus G'']-[M,\widetilde{G}],$$
where $\widetilde{G}$ is obtained from $G$ by recoloring all
components of $G''$ in a new, $(n+1)$-th, color.
The defect of the coloring for each of $G\sminus G'$,
$G\sminus G''$ and $\widetilde{G}$ is less than $d$.
The statement follows by the induction on $d$.
\end{proof}

Thus, without a loss of generality, one can assume
that the colorings in the definition of $\cFY_n(M)$
are simple.

\subsection{The filtration $\cFY_n$}
Further on we will be mainly interested in $M$ being
a $\Z$HS.
Denote $\cFY_n(S^3)$ by $\cFY_n$.

\begin{lem}\label{lem_sphere}
Let $G$  be an $n$--colored Y--link in a $\Z${\rm HS} $M$.
Then $[M,G]\in\cFY_n$.
\end{lem}

\begin{proof}
It suffices to prove, that there exist two Y--links
$G'$ and $G''$ in $S^3$, each colored in at least
$n$ colors, such that $[M,G]=[S^3,G'']-[S^3,G']$.

Since $M$ is a $\Z$HS, it may be obtained from $S^3$
by Y--surgery on a Y--link $G$ in $S^3$.
Let $G'$ be an $n$--colored Y--link in $S^3\sminus G$,
such that its image under Y--surgery on $G$ is
isotopic to $G$.
Let $G''$ be the $(n+1)$--colored Y--link, obtained
from $G'$ by an addition of $G$, with all components
of $G$ colored in a new color.
Then $[S^3,G'\cup G]=[S^3_G,G]+[S^3,G']$ and hence
$[S^3_G,G]=[S^3,G'']-[S^3,G']$.
\end{proof}

\subsection{The filtrations $\cFa_n$ and $\cFb_n$}
A link $L$ in a $\Z$HS is {\em algebraically split\/} if
all pairwise linking numbers of its components vanish.
Such a link is {\em boundary} if all its components
bound non-intersecting surfaces.
A framing of $L$ is {\em unimodular} if the self-linking
of each component is $\pm1$.
Surgery on (any sublink of) a unimodular algebraically
split link gives again a $\Z$HS.
Using these classes of links, T~Ohtsuki \cite{Oh}
and the first author \cite{Ga} introduced two different
filtrations on a vector space generated by all $\Z$HS.
Below we describe the corresponding filtrations on
the free abelian group generated by all $\Z$HS.

For a framed link $L$ in $M$, denote by $M_L$ the result
of surgery of $M$ along $L$.
In $\cF$, let $$[M,L]=\sum_{L'\subset L}
(-1)^{|L'|}M_{L'},$$ where $|L'|$ is the number of
components of a sublink $L'$.
Let $\cFa_n$ (respectively $\cFb_n$) be the subgroup of
$\cF$, generated by $[M,L]$ for all  algebraically split
(respectively boundary) unimodular links $L$ in
$\Z$HS $M$ with at least $n$ components.
As shown in \cite{Oh} and \cite{Ga}, in the definitions of
these filtrations it suffices to consider only $M=S^3$.

\subsection{The filtration $\cFbl_n$}
In \cite{GL}, J~Levine and the first author considered
another filtration, based on a notion of a blink.
A {\em blink} is a framed link $B$ with a given splitting of
the set of its components into pairs
$(B_1^-,B_1^+),\dots,(B_n^-,B_n^+)$, such that:
\begin{itemize}
\item each pair $(B_i^-,B_i^+)$ bounds an orientable surface
$\GS_i$, so that $\GS_i\cap\GS_j=\emptyset$ for $i\ne j$;
\item the surface $\GS_i$ induces a preferred framing of
$B_i^-$ and $B_i^+$; this framing should differ by $\pm1$
from the given framing of $B_i^\pm$.
\end{itemize}
A {\em sub-blink} $B'$ of $B$ is obtained from $B$ by a
removal of some pairs $(B_i^-,B_i^+)$.
In $\cF$, put $$[M,B]=\sum_{B'\subset B} (-1)^{|B'|}M_{B'},$$
where $|B'|$ is the number of pairs in a sub-blink $B'$.
Let $\cFbl_n(M)$ be the subgroup of $\cF$, generated by
$[M,L]$ for all blinks $L$ in a manifold $M$ with at
least $n$ pairs of components.

In the same paper \cite{GL} several other filtrations were
introduced, using different subgroups of the mapping class
group. Each of them was shown to be equivalent to one of the
above filtrations.

\subsection{The main results}\label{sub_main}
We describe a topological calculus of surgery
on Y--links, developed by the second author.
Using this calculus, we describe the structure
of the graded groups $\cFY_n/\cFY_{n+1}\oZ$,
where $\Z[1/2]=\{n/2^m|n,m\in\Z\}$ is the group
of binary rational numbers.

We also compare different filtrations:
\begin{itemize}
\item $\cFY_n(M)=\cFbl_n(M)$;
\item $\cFY_{2n}\subset\cFa_{3n}$ and
$\cFY_{2n}\subset\cFb_{n}$;
\item $\cFY_{2n-1}\oZ=\cFY_{2n}\oZ$;
\item $\cFY_{2n}\oZ=\cFa_{3n}\oZ$ and
$\cFY_{2n}\oZ=\cFb_n\oZ$.
\end{itemize}

Let $F$ be a free abelian group on a set $S$ of
generators, equipped with a decreasing filtration
$F=F_0\supset F_1\supset\dots$. Given an abelian group
$A$, a function $S\to A$ is called a {\em finite type
invariant} of degree at most $d$, if its extension
to a homomorphism $F\to A$ vanishes on $F_{d+1}$.
Each filtration ${\cFa_n}$, ${\cFb_n}$, ${\cFbl_n}$ and
${\cFY_n}$ defines a notion of a finite type invariant
on the set of all $\Z$HS.
The above comparison of filtrations implies the
following theorem announced in \cite[Theorem 1]{Gu2}.
If 2 is invertible in $A$, these definitions are
equivalent, with the following relation of degrees:
$$2d^{as}=6d^b=3d^{bl}=3d^Y.$$

In \cite{GL}, weaker results  $\cFb_n \otimes\BQ
\subset\cFbl_{2n}\otimes\BQ=\cFa_{3n}\otimes\BQ$
were obtained by using different methods that involved
the study of the mapping class group and several of its
subgroups  (building on the results of D Johnson \cite{Jo,Jo2}).

\section{Surgery on clovers}\label{sec_clovers}
In this section we introduce a generalization of Y--links,
which we call clovers. They turn out to be quite useful
from a technical point of view, and are closely related
to uni-trivalent graphs appearing in the study of the
graded quotients $\cFY_n/\cFY_{n+1}$. Similar objects
were called {\em allowable graph claspers}\footnote{See 
{\tt http://www.dictionary.com} for some unexpected
meanings of the word clasper.}
%{\bf clasper} {\em n.} \begin{enumerate}
%\item One that clasps.
%\item Any of the appendages of the male of certain
%insects and crustaceans that are used during copulation
%to hold the female. \item A posterior extension on each
%of the pelvic fins of male elasmobranch fishes that
%aids in the transmission of sperm during copulation.
%\end{enumerate} (Dictionary of the English language)}
by K Habiro \cite{Ha}.
We start by recalling some basic facts about surgery
(see \cite{K}).
\subsection{Surgery on links}
\begin{figure}[htb]
\centerfig{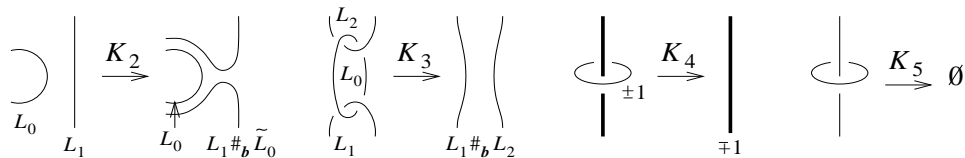,width=5in}
\caption{Kirby moves}
\label{fig_kirby}
\end{figure}

\begin{figure}[htb]
\centerfig{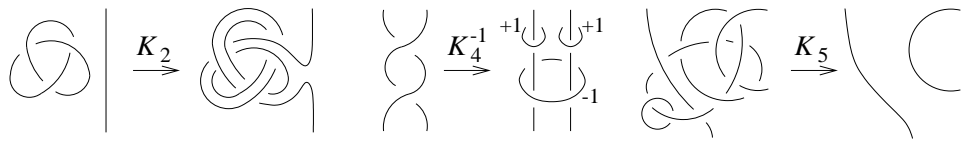,width=5in}
\caption{Examples of Kirby moves}
\label{fig_kirby_ex}
\end{figure}

We will call framed links $L$ and $L'$ in a manifold $M$
{\em surgery equivalent}, and denote $L\sim L'$, if $M_L$
is diffeomorphic to $M_{L'}$.
Recall (see \cite{K}), that links $L$ and $L'$ in $S^3$
are surgery equivalent, iff one can pass from $L$ to $L'$
by a sequence of {\em Kirby moves} $K_1$ and $K_2$:
\begin{itemize}
\item[$K_1$:]
Add to $L$ a small $\pm 1$ framed unknot, unlinked with
the other components of $L$.
\item[$K_2$:]
Add $L_0$ to $L_1$ by replacing $L_1$ with $L_1\#_b\thinspace
\widetilde{L_0}$, where $\#_b$ is a band connected sum and
$\widetilde{L_0}$ is a push-off of $L_0$ along the framing.
\end{itemize}
It is convenient to introduce three additional moves
$K_3$--$K_5$, which can be expressed via $K_1^{\pm 1}$
and $K_2$ (see \cite{K} for $K_4$, $K_5$, and \cite{MP}
for $K_3$):
\begin{itemize}
\item[$K_3$:]
Let $L_0$ be a closed, $0$--framed component of $L$ bounding
an embedded disk $D$, which intersects $L\sminus L_0$ in
exactly two points, belonging to different components $L_1$
and $L_2$ of $L$.
Replace $L_0 \cup L_1 \cup L_2$ by
$L_1\#_b\thinspace L_2$, where the band $b$ intersects $D$
along the middle line of $b$.
\item[$K_4$:] Delete a $\pm 1$--framed unknot, at the expense
of the full left- or right-hand twist on the strings linked
with it.
\item[$K_5$:]
Let $L_0$ be a closed, $0$--framed component of $L$ which bounds
an embedded disk $D$. Suppose that $D\cap (L\sminus L_0)$
consists of exactly one point lying on some component
$L_1\subset L$.

Delete $L_0$ and $L_1$ from $L$.

\end{itemize}
The moves $K_2$--$K_5$ are shown in Figure \ref{fig_kirby}.
Here, and throughout the paper, we think about a surgery
presentation of $M$ and use thick solid lines to depict an
arbitrary union of surgery components comprising $M$ (possibly
together with some edges of embedded graphs).
Also, we use a usual convention that all depicted links, or
graphs, coincide outside a ball shown in the picture.
In all figures the framings are assumed to be orthogonal to
the plane of the picture, unless explicitly indicated
otherwise.
Some examples of the moves $K_2$, $K_4$ and $K_5$ are given
in Figure \ref{fig_kirby_ex}.
By a repeated application of $K_3$ we obtain the following
relation (which explains the term {\em Borromean surgery},
used in \cite{Ma} for Y--surgery in a manifold $M$, and
the term {\em $\GD$--move}, used in \cite{MN}):
\begin{lem}\label{lem_borro}
$$\fig{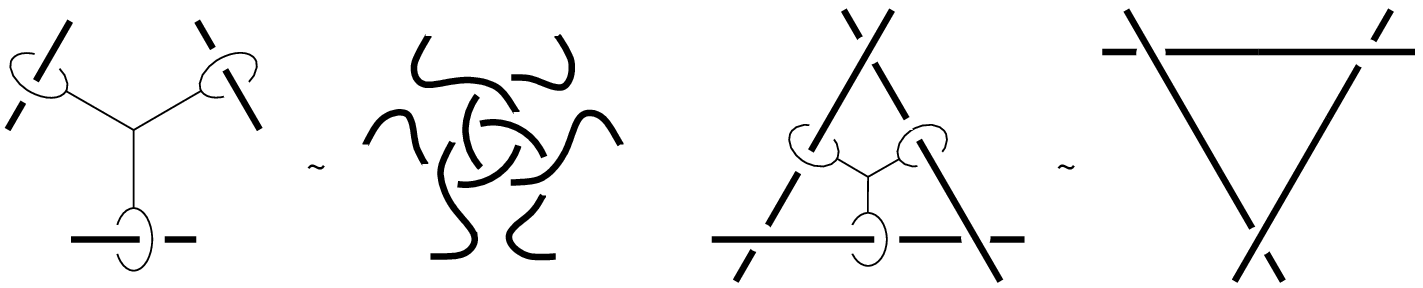,height=0.9in}$$
\end{lem}

\subsection{Clovers}
Let $G$ be a framed trivalent graph, smoothly embedded
into a manifold $M$.
We will call $G$ a {\em clover}, if:
\begin{itemize}
\item Each point of $G$ has a neighborhood,
diffeomorphic to a neighborhood of some point of a
standard Y--graph (including the framing), by an
orientation-preserving diffeomorphism.
\item Different edges meeting in a vertex have different
tangent lines.
\end{itemize}
If both endpoints of an edge coincide, the two
corresponding tangent lines may coincide.
In this case the edge is called a {\em leaf}, and the
incident vertex is called an {\em external vertex}.
A vertex not incident to a leaf is {\em internal}.
We call an edge {\em internal}, if both its endpoints
are internal vertices. A {\em degree} of $G$ is the
number of internal vertices.
For example, a Y--graph is a degree 1 clover.
It has one internal vertex, three external vertices
and three leaves, see Figure \ref{fig1}a.

Clovers of degree 0 play an important role in an
analogous theory of finite type invariants of links.
However, we will exclude them from the theory of
finite type invariants of 3--manifolds.
Therefore all throughout the paper we will always
assume that each connected component of a clover
is of degree at least one.

In figures we will always assume that the framing of
$G$ is orthogonal to the plane of the picture.

\subsection{Surgery on clovers}
\begin{figure}[htb]
\centerfig{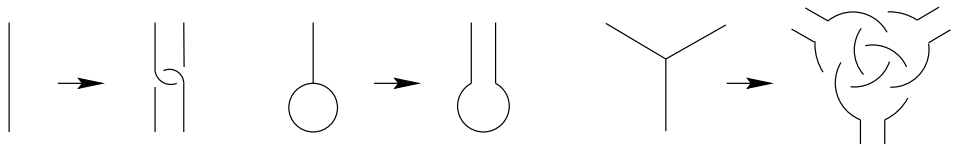,height=0.75in}
\caption{Construction of the surgery link from a clover}
\label{fig2}
\end{figure}

\begin{figure}[htb]
%\vskip1in
\centerfig{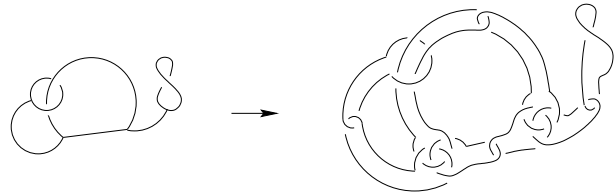,height=1.0in}
\caption{A clover and the corresponding link}
\label{fig3}
\end{figure}

Let $G$ be a clover in $M$. We construct a framed
link $L(G)$ in a small neighborhood of $G$. In Figure
\ref{fig2} we show the structure of $L(G)$ near an edge, a
leaf, and an internal vertex of $G$. The framings of the
fragments of $G$ and $L(G)$ appearing in Figure \ref{fig2}
are assumed to be orthogonal to the plane of the picture.
The construction of $L(G)$ is illustrated in Figure
\ref{fig3}.
The surgery of $M$ along $L(G)$ is called a {\em surgery
along $G$}, and the manifold $M_{L(G)}$ is denoted by $M_G$.
We call two clovers $G$ and $G'$ in $M$ {\em surgery
equivalent}, and denote $G\sim G'$, if $M_G$ is
diffeomorphic to $M_{G'}$.

Denote by $G_1,\dots,G_n$ the connected components of $G$.
For an $n$--tuple $\Gs=\{\Gs_1,\dots,\Gs_n\}\in\{0,1\}^n$
we put $$[M,G]=\sum_\Gs (-1)^{|\Gs|}M_{G(\Gs)},$$ where
$G(\Gs)=\bigcup_{i:\Gs_i=1}G_i$ and $|\Gs|=\sum_i\Gs_i$.
If the degree of each $G_i$ is 1, we recover the
definition of $[M,G]$ for Y--links.

Suppose that a leaf or a cycle of edges of $G$ in $M$
bounds an embedded surface $\GS$ in $M$.
Then the intersection number of its push-off along the
framing with $\GS$ does not depend on the choice of
$\GS$ and (by a slight abuse of terminology) will be
also called the framing of this leaf or cycle.
A leaf is {\em trivial}, if it is 0--framed and
bounds a disc $D$ in $M$ whose interior does not
intersect $G$.
\begin{lem}\label{lem_disc}
Let $G\subset M$ be a connected clover, which
contains a trivial leaf $l$.
Then surgery on $G$ preserves a neighborhood of
$G\cup D$, where $D$ is the disc bounding $l$.
In particular, $M_G=M$ and $[M,G]=0$:
$$\fig{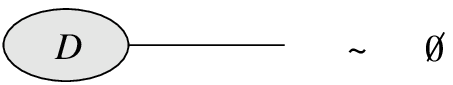,width=1.8in}$$
\end{lem}
\begin{proof}
Follows from the construction of $L(G)$ by an application
of $K_5$.
\end{proof}

An edge of a clover is a {\em trivial loop}, if both its
endpoints coincide, it is 0--framed, and bounds a disc $D$
whose interior does not intersect $G$.
\begin{lem}\label{lem_loop}
Let $G\subset M$ be a connected clover, which contains
a trivial looped edge $e$.
Then surgery on $G$ preserves a neighborhood of
$G\cup D$, where $D$ is the disc bounding $e$.
In particular, $M_G=M$ and $[M,G]=0$.
\end{lem}
\begin{proof}
Follows from the construction of $L(G)$ by an application
of $K_5$:
$$\fig{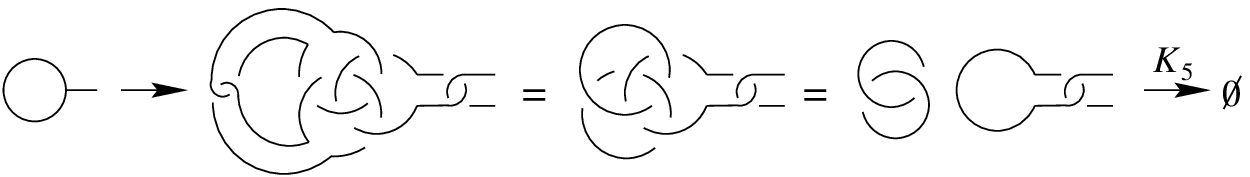,height=0.7in}$$
\end{proof}
\begin{thm}\label{thm_edgecut}
Let $G$ be a clover in $M$, and let $G'$ be
obtained from $G$ by cutting an internal edge of $G$ and
inserting there two small Hopf-linked leaves:
$$\centerfig{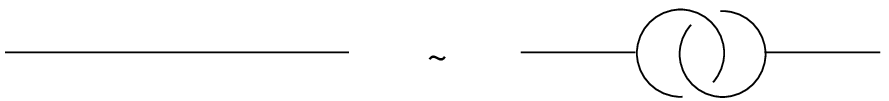,width=3.5in}$$
Then $G'\sim G$ and $[M,G']=\pm[M,G]$, where the sign is
negative if the edge is splitting $G$, and positive otherwise.
\end{thm}
\begin{proof}
The equality $G\sim G'$ follows from the construction
of $L(G)$ by an application of $K_3$.
Topologically, this corresponds to cutting the corresponding
solid handlebody introducing a pair of complimentary handles
as shown below:
$$\centerfig{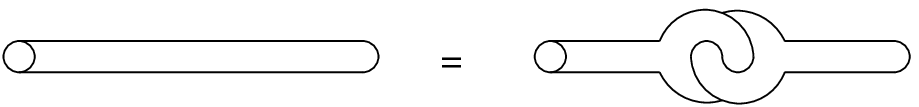,width=3.5in}$$
Let $G_0$ be the component of $G$ containing this
edge, and $G_1$, $G_2$ be the components of $G'$
replacing $G_0$.
Note that each subclover of $G'$ which contains exactly
one of the components $G_1$ and $G_2$, has a trivial leaf.
Thus by the definition of $[\,\cdot\,,\,\cdot\,]$ and
Lemma \ref{lem_disc} $[M,G']=[M,G]$ if $G'$ has the same
number of components as $G$, and $[M,G']=-[M,G]$ otherwise.
\end{proof}
\begin{cor}\label{cor_graph}
Let $G$ be a degree $n$ clover.
Then $[M,G]\in\cFY_n(M)$.
\end{cor}

\section{Topological calculus of clovers}
\label{sec_calculus}
Below we describe some important moves $Y_1$--$Y_4$ on
Y--links, which are shown to preserve the surgery
equivalence classes. One can apply these moves also to
clovers, in view of Theorem \ref{thm_edgecut}.

\begin{figure}[htb]
\centerfig{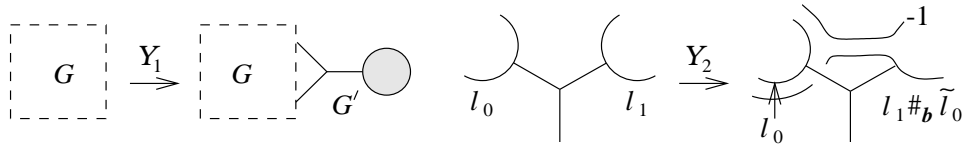,width=5.0in}
\caption{Blow-up and leaf slide}
\label{fig_Ymoves12}
\end{figure}

\begin{figure}[htb]
\centerfig{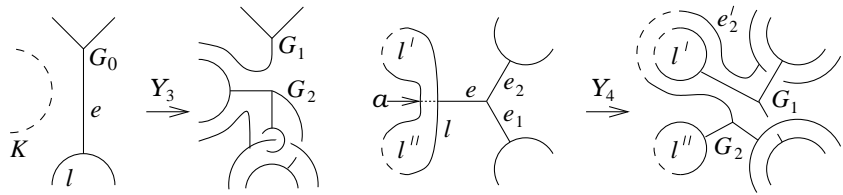,height=1.0in}
\caption{Edge slide and leaf-cutting}
\label{fig_Ymoves34}
\end{figure}

Moves $Y_1$ and $Y_2$ (similar to the Kirby moves $K_1$,
$K_2$) are shown in Figure \ref{fig_Ymoves12}:
\begin{itemize}
\item[$Y_1$:]
Add to $G$ a connected clover $G'$, which has a 0--framed
leaf $l$ bounding an embedded disc whose interior does
not intersect $G\cup G'$.
\item[$Y_2$:]
Add a leaf $l_0$ to another leaf $l_1$ of the same Y--graph
along a band $b$ by replacing $l_1$ with
$l_1\#_b\thinspace\widetilde{l_0}$, where $\widetilde{l_0}$
is a push-off of $l_0$ along the framing; then change the
framing by $-1$.
\end{itemize}
Let $l$ be a leaf of a Y--graph $G_0$ and $e$ be the
adjacent edge. Let $K$ be a knot.
The move $Y_3$ is sliding an edge $e$ along $K$, see
Figure \ref{fig_Ymoves34}:
\begin{itemize}
\item[$Y_3$:]
Denote by $G_1$ Y--graph obtained from $G_0$ by adding $K$
to $e$ along a band $b$. Construct a Y--graph $G_2$ as
shown in Figure \ref{fig_Ymoves34}: one of its leaves is
a push-off of $K$ (and the adjacent edge goes along $b$),
another is a push-off of $l$, and the third leaf is a
0--framed unknot linked once with $l$.
Replace $G_0$ by $G_1\cup G_2$.
\end{itemize}
Let $l$ be a leaf of a Y--graph $G_0$, $e$, $e_1$ and
$e_2$ be its edges, with $e$ adjacent to $l$.
An arc $a$, starting at $e\cap l$ and ending on $l$,
cuts $l$ into two parts $l'$ and $l''$.
The move $Y_4$ is cutting the leaf $l$ along $a$, see
Figure \ref{fig_Ymoves34}:
\begin{itemize}
\item[$Y_4$:]
Let $G_1$ be an Y--graph obtained from $G_0$ by
replacing the leaf $l$ by $l'\cup a$.
Let $G_2$ be obtained from $G_0$ by replacing the
leaf $l$ by $l''\cup a$, adding $l'$ to $e_2$ along
the band $e$, and taking a push-off copy $e'_2$, as
shown in Figure \ref{fig_Ymoves34}.
Replace $G_0$ by $G_1\cup G_2$.
\end{itemize}

A convenient way to draw the moves $Y_1$--$Y_4$ is by
depicting them in a standard handlebody, which then
can be embedded into a manifold in an arbitrary way.
Below is such an interpretation of the moves $Y_3$,
$Y_4$ as moves in a standard handlebody of genus 4:
\begin{figure}[htb]
\centerfig{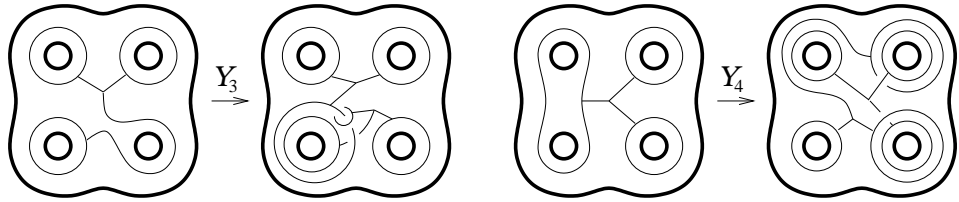,width=5in}
\caption{$Y_3$ and $Y_4$ as moves in a handlebody}
\label{fig_handle34}
\end{figure}

\begin{figure}[htb]
\centerfig{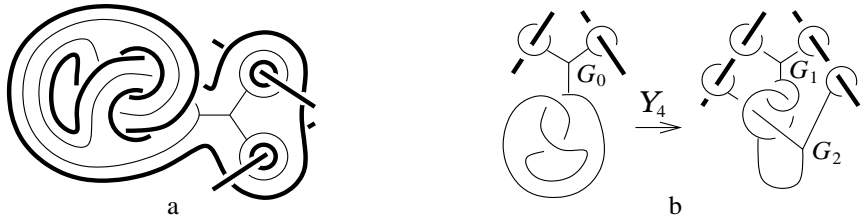,width=4.5in}
\caption{An embedding and the induced $Y_4$ move}
\label{fig_handle}
\end{figure}

For example, when the handlebody is embedded so that
one of the handles links another handle as in Figure
\ref{fig_handle}a, the corresponding move $Y_4$ is
shown in Figure \ref{fig_handle}b.
\begin{thm}\label{thm_equiv}
The moves $Y_1$--$Y_4$ preserve the classes of surgery
equivalence of Y--links in a manifold $M$.
\end{thm}
\begin{proof}
It suffices to prove the surgery equivalence of Y--links
obtained by the moves $Y_1$--$Y_4$ in the standard
handlebody. Instead of drawing the handlebodies we
will draw thick lines passing through the handles
(encoding a set of surgery and clover components),
similarly to Lemma \ref{lem_borro}.

By Lemma \ref{lem_disc}, $Y_1$ is surgery equivalence.
To verify $Y_2$, depict the Borromean linking with
one component passing on the boundary of an embedded
surface with two handles, and then twist one of the
handles along the other:
$$\fig{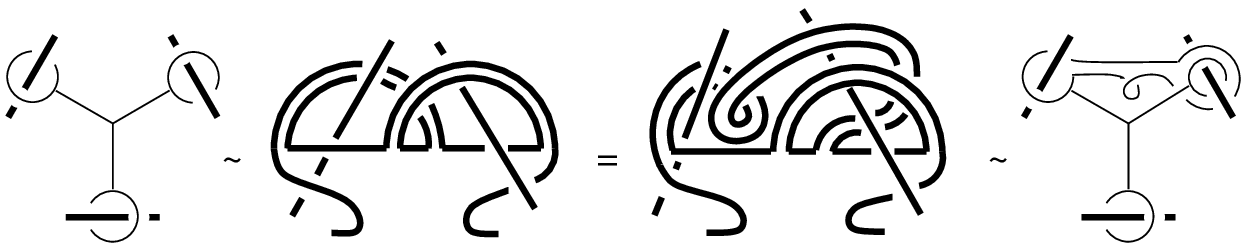,height=1.0in}$$
To verify $Y_3$, use an isotopy and Lemma \ref{lem_borro}:
$$\fig{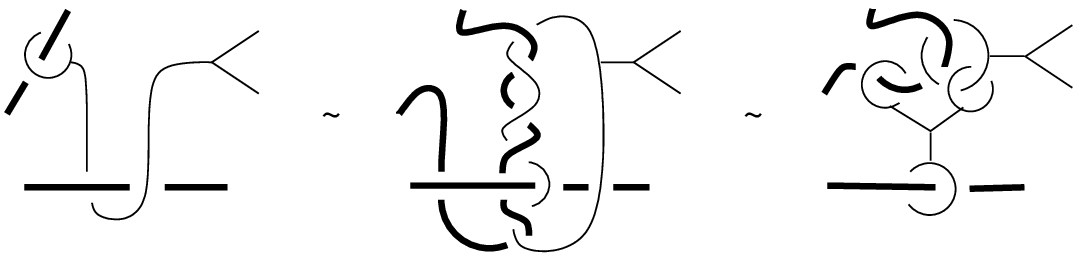,height=1.0in}$$
The verification of $Y_4$ is similar:
$$\fig{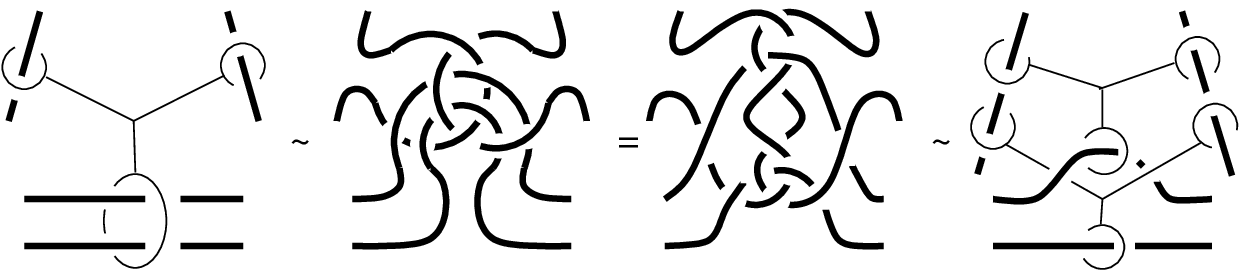,height=1.0in}$$
\end{proof}

\begin{figure}[htb]
\centerfig{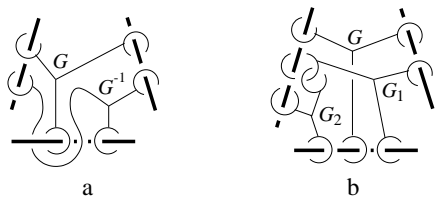,height=1.0in}
\caption{Two presentations of the inverse of a Y--graph}
\label{fig_inverse}
\end{figure}

The following theorem was announced in \cite{Gu2} in
case of a Y--graph $G$.
\begin{thm}\label{thm_inverse}
For any clover $G$ in a manifold $M$, there exists a
clover $G^{-1}$ in a neighborhood of $G$, such that
the result of surgery on both clovers is the original
manifold: $M_{G\cup G^{-1}}=M$.
The construction of $G^{-1}$ for a Y--graph $G$ is shown
in Figure \ref{fig_inverse}a. Another presentation of
$G^{-1}$ by a 2--component Y--link $G_1\cup G_2$ is shown
in Figure \ref{fig_inverse}b.
\end{thm}
\begin{proof}
Let us first verify the statement for a Y--graph.
By an isotopy and a subsequent application of $Y_4$
and $Y_1$, we get
$$\fig{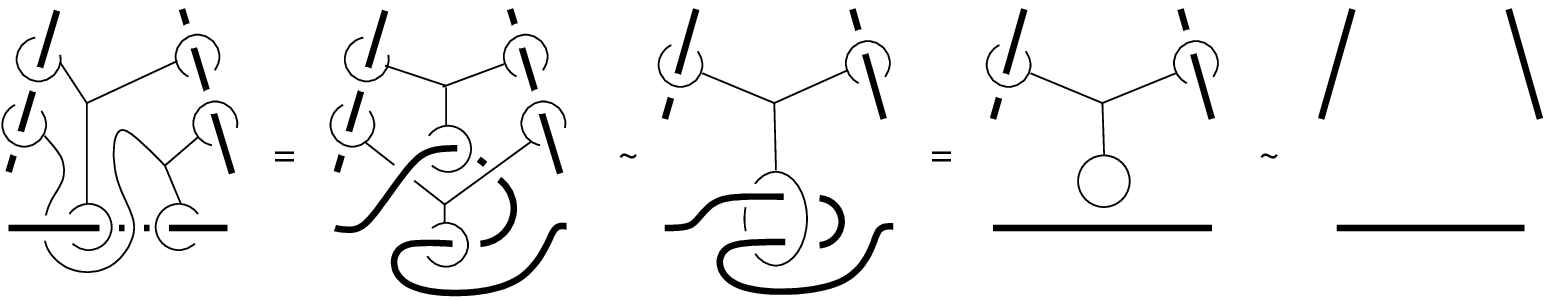,width=5.0in}$$
The general case now follows by Theorem \ref{thm_edgecut}.
Finally, $G^{-1}\sim G_1\cup G_2$ by $Y_3$.
\end{proof}

\section{The structure of the graded spaces $\cG_n$}
\label{sec_graded}
This section is devoted to the study of the graded
spaces$$\cG_n(M)=\cFY_n(M)/\cFY_{n+1}(M).$$
Denote $[M,G]\eq{n}[M,G']$ iff
$[M,G]-[M,G']\in\cFY_{n+1}(M)$.

\subsection{Graded versions of $Y_2$--$Y_4$}
Using \ref{thm_edgecut} and Lemma \ref{lem_disc}, we
obtain the following immediate corollaries of Theorem
\ref{thm_equiv}.

\begin{cor}\label{cor_Y2}
Let $G$ be a clover of degree $n$ in a manifold $M$,
and let $G'$ be obtained from $G$ by sliding one of
its leaves along an adjacent leaf by $Y_2$.
Then $[M,G']=[M,G]$, and hence $[M,G']\eq{n}[M,G]$.
\end{cor}
Let $G_0$ be a Y--graph in $M$, and let $G_1\cup G_2$
be obtained from $G_0$ by an application of $Y_3$, see
Figure \ref{fig_Ymoves34}.
Note that $M_{G_2}=M$ by $Y_1$, thus
$$[M,G_0]-[M,G_1]=M_{G_1\cup G_2}-M_{G_1}=
[M,G_1\cup G_2]\in\cFY_2(M).$$
Hence, using Theorem \ref{thm_edgecut} we obtain:
\begin{cor}\label{cor_Y3}
Let $G$ be a clover of degree $n$ in a
manifold $M$, and let $K\subset M$ be a knot.
Let $G'$ be obtained from $G$ by sliding an edge
of $G$ along $K$.
Then $[M,G']\eq{n}[M,G]$.
\end{cor}
In Lemma \ref{lem_slide} we will strengthen this
result by computing the difference $[M,G']-[M,G]$
in $\cG_{n+1}(M)$.

\begin{figure}[htb]
\centerfig{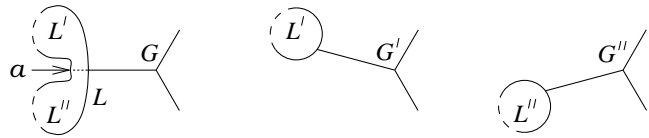,height=0.7in}
\caption{Splitting a leaf}
\label{fig_split}
\end{figure}

Let $G_0$ be a Y--graph in $M$, and let $G_1\cup G_2$
be obtained from $G_0$ by an application of $Y_4$,
see Figure \ref{fig_Ymoves34}.
By the previous Lemma, the value of $[M,G_2]$ in $\cG_1(M)$
does not change when we slide its edge $e'_2$ along a knot.
Thus we can pull $e'_2$ off $l'$, making it back into $e_2$.
Hence, using again Theorem \ref{thm_edgecut} we obtain:
\begin{cor}\label{cor_Y4}
Let $G$ be a clover of degree $n$ in a
manifold $M$ and $l$ be a leaf of $G$.
An arc $a$ starting in the external vertex incident
to $l$ and ending in other point of $l$, splits $l$
into two arcs $l'$ and $l''$.
Denote by $G'$ and $G''$ the clovers obtained from $G$
by replacing the leaf $l$ with $l'\cup a$ and
$l''\cup a$ respectively, see Figure \ref{fig_split}.
Then $[M,G]\eq{n}[M,G']+[M,G'']$.
\end{cor}
\subsection{The dependence on framings}
Theorem \ref{thm_inverse} allows us to deduce
the dependence of $[M,G]$ on the framings of edges.
\begin{lem}\label{lem_frame}
Let $G$ be a clover of degree $n$ in a manifold $M$.
Let $G'$ be obtained from $G$ by twisting the framing
of an edge by a half twist.
Then $[M,G']\eq{n}-[M,G]$.
\end{lem}
\begin{proof}
Let $G$ and $G_1\cup G_2\sim G^{-1}$ be the Y--graphs
depicted in Figure \ref{fig_inverse}b.
Note that $G_1$ looks exactly like $G$, except for
the way its lower leaf links the thick line.
Turning this leaf to the same position changes the
framing of the adjacent edge by a half twist:
$$\fig{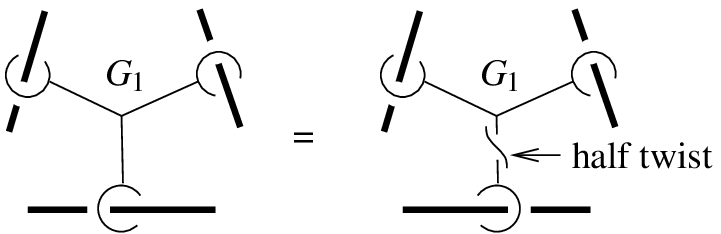,height=0.85in}$$
By Theorem \ref{thm_inverse}, $M_{G\cup G_1\cup G_2}=M$;
also, $M_{G\cup G_2}=M_G$ and $M_{G_2}=M$ by Lemma
\ref{lem_disc}. Thus
\begin{eqnarray}
[M,G]+[M,G_1]&=&2M-M_G-M_{G_1}\nonumber\\&=&
[M,G\cup G_1]+[M,G_1\cup G_2]-[M,G\cup G_1\cup G_2]\eq{1}0\nonumber
\end{eqnarray}
and the lemma follows.

Alternatively, one can show that $[M,G]+[M,G^{-1}]\eq{1}0$
for a Y--graph $G^{-1}$ depicted in Figure \ref{fig_inverse}a
and pull its edge off the lower thick line by \ref{cor_Y3}
to obtain once again $G_1$.
\end{proof}
\begin{cor}\label{cor_frame}
Let $G$ be a clover of degree $n$ in $M$.
Adding a kink, ie, a full twist to the framing of an
edge preserves an $n$--equivalence class of $[M,G]$.
\end{cor}

\begin{figure}[htb]
\centerfig{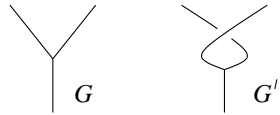,height=0.6in}
\caption{The AS relation}
\label{fig_order}
\end{figure}

Now, note that the following two Y--graphs are isotopic:
$$\fig{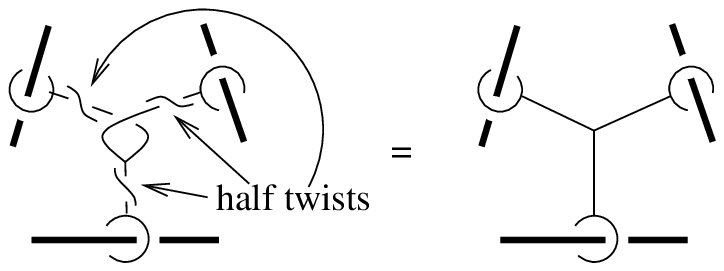,height=1.2in}$$
Thus, applying Lemma \ref{lem_frame} and Theorem
\ref{thm_edgecut}, we obtain:
\begin{cor}\label{cor_order}
Let $G$ and $G'$ be clovers of degree $n$
in a manifold $M$, which coincide everywhere except
for a fragment shown in Figure \ref{fig_order}.
Then $[M,G']\eq{n}-[M,G]$.
\end{cor}

\subsection{Simplifying the leaves}\label{sub_simple}
Denote $\cG_n(S^3)$ by $\cG_n$.
We want to show that the space $\cG_n\oZ$ is generated
by clovers with only internal vertices, ie, without
leaves.

A leaf of a clover $G$ is {\em special}, if it is
$\pm1$--framed and bounds an embedded disc, whose
interior does not intersect $G$.
A leaf $l$ of a clover $G$ is {\em simple}, if it either
bounds an embedded $i$--framed disc $D$ whose interior
intersects $G$ in at most one point, or is special.

For $M=S^3$, the graded space $\cG_n$ is generated by
clovers of degree $n$ all leaves of which are simple.
Indeed, suppose we are given an arbitrary clover of
degree $n$ in $S^3$.
We may split each of its non-trivial leaves into small
pieces and apply Corollary \ref{cor_Y4} to present $G$
in $\cG_n$ as a linear combination of clovers with
simple leaves.

\begin{figure}[htb]
\centerfig{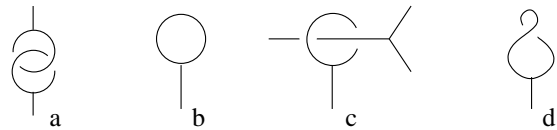,height=0.65in}
\caption{Four types of simple leaves}
\label{fig_simple}
\end{figure}

Suppose that all leaves of clover $G$ in a manifold $M$
are simple.
If two simple leaves are linked with each other (see
Figure \ref{fig_simple}a), we replace them by a new
edge by Theorem \ref{thm_edgecut}.
Suppose that at least one simple leaf $l$ still
remains after this procedure.
We will show below that in this case
$[M,G]\in\cFY_{n+1}(M)\oZ$.

We should consider three cases, see Figure
\ref{fig_simple}b--d.
Firstly, it may be that $l$ is trivial.
Secondly, it may be that $l$ is 0--framed, and the
disc $D$ intersects an edge of $G$.
Thirdly, it may be that $l$ is special.
In the first case, $[M,G]=0$ by Lemma \ref{lem_disc}.
The second case can be reduced to the first case after
unlinking the edge from the corresponding leaf by
Corollary \ref{cor_Y3} (sliding the edge along a small
unknot linked once with the leaf).
Thus we readily obtain:
\begin{lem}\label{lem_fragm1}
Let $G$ be a clover of degree $n$ in a manifold $M$.
Suppose that $G$ contains a 0--framed leaf bounding
an embedded disc, whose interior intersects $G$ in
exactly one point, belonging to an edge.
Then $[M,G]\eq{n}0$.
\end{lem}

In the third case we encounter 2--torsion.

\begin{lem}\label{lem_fragm2}
Let $G$ be a clover of degree $n$ in a manifold $M$,
which contains a special leaf.
Then $2[M,G]\eq{n}0$.
\end{lem}
\begin{proof}
Rotating the special leaf we can change the framing
of the adjacent edge while preserving the isotopy
class of $G$. Thus by Lemma \ref{lem_frame}
$[M,G]\eq{n}-[M,G]$ and the lemma follows.
\end{proof}

Over the integers, we have the following inclusion:
\begin{lem}\label{lem_fragm2'}
Let $G$ and $l$ be as in Lemma \ref{lem_fragm2}.
Suppose that the connected component of $G$ which
contains $l$ is of degree at least two.
Then $[M,G]\eq{n}0$.
\end{lem}
\begin{proof}
Cutting the neighboring internal edge we obtain
two new leaves $l_1$ and $l_2$.
Without a loss of generality suppose that the
framing of $l$ is $+1$.
Then slide $l$ along $l_1$ by $Y_2$ as shown below:
$$\fig{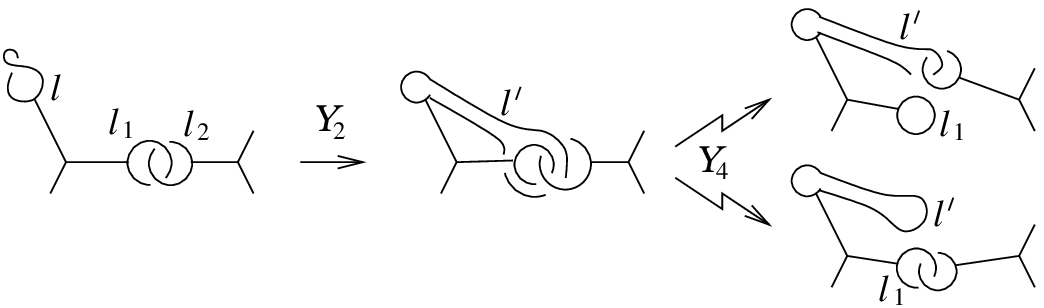,height=1.2in}$$
Notice that $Y_2$ changes the framing of a leaf
by $-1$, so the new leaf $l'$ is 0--framed.
Splitting $l_2$ as shown above by Corollary \ref{cor_Y4},
we obtain two clovers each of which contains a
trivial leaf (either $l'$ or $l_1$), and so can be
removed by $Y_1$.
\end{proof}

\subsection{The IHX relation}\label{sec_IHX}
\begin{figure}[htb]
\centerfig{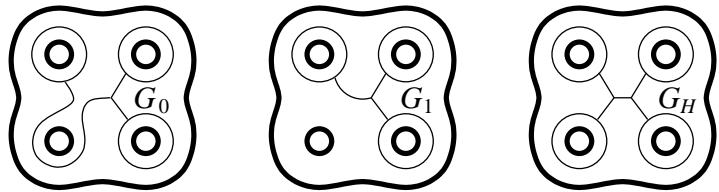,height=1.0in}
\caption{Sliding an edge and computing the difference}
\label{fig_slide1}
\end{figure}

\begin{figure}[htb]
\centerfig{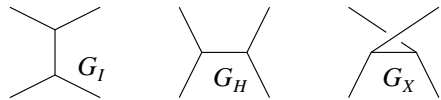,height=0.5in}
\caption{The IHX relation}
\label{fig_ihx}
\end{figure}

Let $G=\GG\cup G_0$ be an $n$--component Y--link which
contains a Y--graph $G_0$, and let $K$ be a knot in $M$.
Choose a band connecting an edge of $G_0$ to $K$.
Sliding this edge of $G_0$ along $K$ by $Y_3$, we obtain
a Y--graph $G_1$, as shown in Figure \ref{fig_Ymoves34}.
Denote $G'=\GG\cup G_1$. By Corollary \ref{cor_Y3},
$[M,G]-[M,G']\in\cFY_{n+1}(M)$.
There is a simple expression for this difference modulo
$\cFY_{n+2}(M)$.
It is easier to visualize the picture in a neighborhood
$N$ of $G_0\cup b\cup K$, which is a genus 4 handlebody
embedded into $M$.
\begin{lem}\label{lem_slide}
Let $G_0$, $G_1$ and $G_H$ be the clovers of Figure
\ref{fig_slide1} in a handlebody $V$  embedded into $M$.
Let $\GG$ be a degree $n-1$ clover in the
complement of $N$.
Put $G=\GG\cup G_0$, $G'=\GG\cup G_1$ and
$G''=\GG\cup G_H$. Then $[M,G]-[M,G']\eq{n+1}[M,G'']$.
\end{lem}
Here is another graphical expression for the above
graphs:
$$\fig{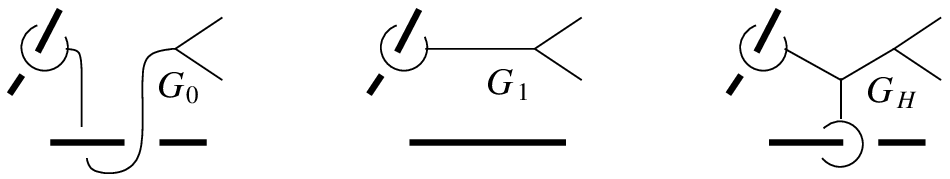,height=0.75in}$$
\begin{proof}
We use $Y_3$ to pass from $G_0$ to $G_1$ and $G_2$
with $G_0\sim G_1\cup G_2$ as in the proof of Theorem
\ref{thm_equiv}.
Then we split the leaf $l$ of $G_1$ by Corollary
\ref{cor_Y4} introducing new Y--graphs $G_1'$ and $G_1''$
with $[M,\GG\cup G_1\cup G_2]\eq{n+1}
[M,\GG\cup G_1'\cup G_2]+[M,\GG\cup G_1''\cup G_2]$:
$$\fig{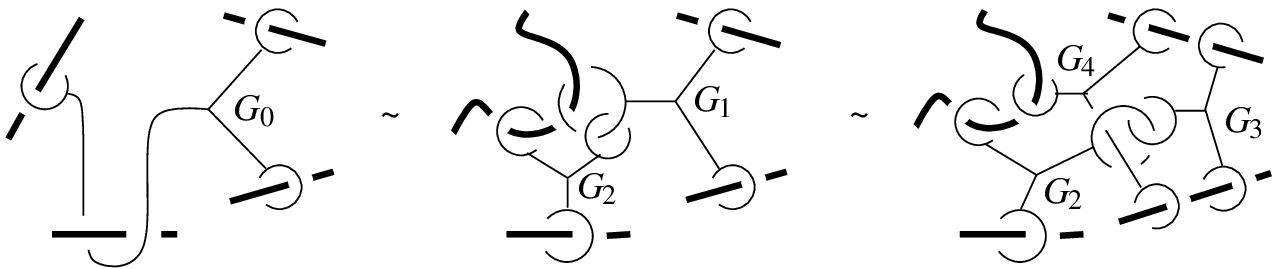,height=1.0in}$$
But $[M,\GG\cup G_1''\cup G_2]=-[M,G'']$ by Theorem
\ref{thm_edgecut}, and $[M,\GG\cup G_1'\cup G_2]=0$
by Lemma \ref{lem_disc}. Hence
$[M,\GG\cup G_1\cup G_2]\eq{n+1}-[M,G'']$.
On the other hand,
$$[M,\GG\cup G_1\cup G_2]=
%sum_{\GG'\subset \GG}
%-1)^{|\GG'|}(M_{\GG'}-M_{\GG'\cup G_1}-
%_{\GG'\cup G_2}+M_{\GG'\cup G_1\cup G_2})=
\sum_{\GG'\subset \GG}(-1)^{|\GG'|}
(-M_{\GG'\cup G_1}+M_{\GG'\cup G_1\cup G_2})=
[M,G']-[M,G]$$
The comparison of two above expressions for
$[M,\GG\cup G_1\cup G_2]$ proves the theorem.
\end{proof}
\begin{thm}\label{thm_IHX}
Let $G_I$, $G_H$ and $G_X$ be clovers
of degree $n$ in a manifold $M$, which coincide
everywhere except for a fragment shown in Figure
\ref{fig_ihx}.
Then $[M,G_I]+[M,G_X]-[M,G_H]\eq{n}0$.
\end{thm}
\begin{proof}
It suffices to prove the statement for $n=2$; the general
case then follows by Theorem \ref{thm_edgecut} and Lemma
\ref{lem_disc}.
Consider a standard Y--graph $G$ in a handlebody $N$
of genus 3 and attach to $N$ an additional handle $h$.
We are to slide all three edges of $Y$ along $h$ in
a genus 4 handlebody $N\cup h$ as indicated in Figure
\ref{fig_ihxpf1}.
Each time we will use Lemma \ref{lem_slide} to compute
the corresponding change of $[M,G]$.
\begin{figure}[htb]
\centerfig{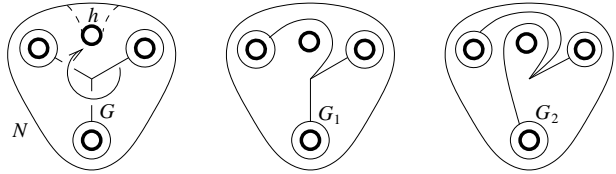,height=1.1in}
\caption{Sliding the edges along a handle}
\label{fig_ihxpf1}
\end{figure}
Sliding the first edge of $G$, we obtain a new Y--graph
$G_1$ with $[M,G_1]-[M,G]\eq{2}[M,G_I]$.
Sliding the next edge of $G$ (or rather of $G_1$), we
obtain a new Y--graph $G_2$ with
$[M,G_2]-[M,G_1]\eq{2}[M,G'_X]$.
After sliding the third edge we return back to $G$ and
get $[M,G]-[M,G_2]\eq{2}[M,G'_H]$.
Here a degree 2 clovers $G_I$, $G'_X$, and $G'_H$ are
shown in Figure \ref{fig_ihxpf2}.
\begin{figure}[htb]
\centerfig{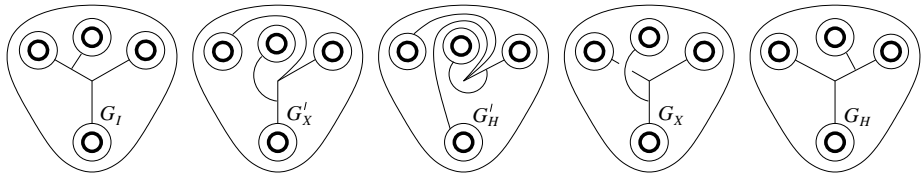,width=5.0in}
\caption{Computing the difference}
\label{fig_ihxpf2}
\end{figure}
Summing up these three equalities we get
\begin{equation}\label{eq_ihx}
[M,G_I]+[M,G'_X]+[M,G'_H]\eq{2}0.
\end{equation}
But $G'_X$ differs from $G_X$ only by an edge slide, and
$G'_H$ differs from $G_H$ by an edge slide and a cyclic
ordering of edges in a vertex, see Figure \ref{fig_ihxpf2}.
Hence $[M,G'_X]\eq{2}[M,G_X]$ by Corollary \ref{cor_Y3},
and $[M,G'_H]\eq{2}-[M,G_H]$ by Corollaries \ref{cor_Y3}
and \ref{cor_order}.
A substitution of two these expressions into
\eqref{eq_ihx} proves the theorem.
\end{proof}

\begin{rem} There is a topological version of the
IHX relation, which may be deduced similarly using
the $Y_3$ and $Y_4$ moves to slide each edge of a
Y--graph through a handle in a handlebody of genus 4.
We leave the details to an interested reader.
\end{rem}

\subsection{Trivalent graphs and $\cG_n$}
\label{sub_graphs}
Consider an abelian group $\tA_k$ freely generated
by abstract (not necessarily connected) trivalent
graphs with $2k$ vertices and without looped edges,
equipped with a cyclic ordering of the incident
edges in each vertex.
%Note that any trivalent graph has an even number of
%vertices, thus $A_{2n-1}=0$ for any $n\in\N$.
Denote by $\cA_k$ the quotient of $\tA_k$ by the
following AS and IHX relations:
\begin{itemize}
\item[AS:]
Let $G'$ be obtained from $G$ by reversing the
cyclic ordering of edges in some vertex, see Figure
\ref{fig_order}.
Then $G'=-G$.
\item[IHX:]
Let $G_I$, $G_H$ and $G_X$ coincide everywhere
except for a fragment shown in Figure \ref{fig_ihx}.
Then $G_I=G_H-G_X$.
\end{itemize}

Denote by $Y$ a Y--graph in $M$ with three special
leaves.
For each graph $G\in\tA_k$ pick an arbitrary
embedding of $G$ into $M$.
Equip it with a framing so that the framing along
each cycle of edges is integer, and take its disjoint
union with $m$ copies of $Y$.
The resulting framed graph in $M$ may be considered
as a clover of degree $n=2k+m$.
Denote it by $\phi(G)$.
Put $\tpsi_n(G)=[M,\phi(G)]$ and extend $\tpsi_n$
to $\tpsi_n\co \oplus_{2k\le n}\tA_k\to\cF$ by linearity.
Note that $\tpsi_n(\tA)\subset\cFY_n$ by Corollary
\ref{cor_graph}.
\begin{thm}\label{thm_graph}
The map $\tpsi_n\co \oplus_{2k\le n}\tA_k\to\cFY_n$
induces a quotient map
$\psi_n\co \oplus_{2k\le n}\cA_k\to\cG_n$,
which is surjective and does not depend on the choice
of $\phi$.
The image of $\oplus\cA_{2k<n}$ is a 2--torsion.
\end{thm}
\begin{proof}
The map $\tpsi_n$ factors through AS and IHX relations
by Corollary \ref{cor_order} and Theorem \ref{thm_IHX}.
The independence on the choice of $\phi(G)$ follows
from Corollary \ref{cor_Y3} and Lemma \ref{lem_frame}.
The surjectivity follows from the results of Section
\ref{sub_simple} and Lemma \ref{lem_loop}.
The torsion result follows from Lemma \ref{lem_fragm2}.
\end{proof}
Note that
$\oplus_{2k\le 2n-1}\cA_k=\oplus_{2k<2n-1}\cA_k$, hence
$\cG_{2n-1}$ consists of a 2--torsion and we obtain:
\begin{cor}
For any $n\in\N$, $\cFY_{2n-1}\oZ=\cFY_{2n}\oZ$.
\end{cor}

\section{The equivalence of the $\cFbl$, $\cFa$
and $\cFY$ filtrations}

\subsection{Equivalence of $\cFbl$ and $\cFY$}
We can present surgery on an arbitrary blink of
genus $g$ as surgery on $g$ blinks of genus one,
slicing the surface into pieces of genus one as
shown below (see also \cite{GL}):
$$\fig{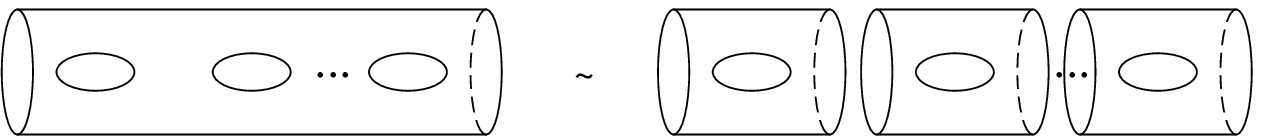,height=0.5in}$$
Also, as noticed already in \cite{Ma}, surgery
on a blink of genus one can be presented as
Y--surgery (and vice versa).
Indeed, depicting the Borromean linking with
one component passing on the boundary of an
embedded surface with two bands, and then
using $K_3^{-1}$, we obtain:
\begin{lem}\label{lem_blink}
$$\fig{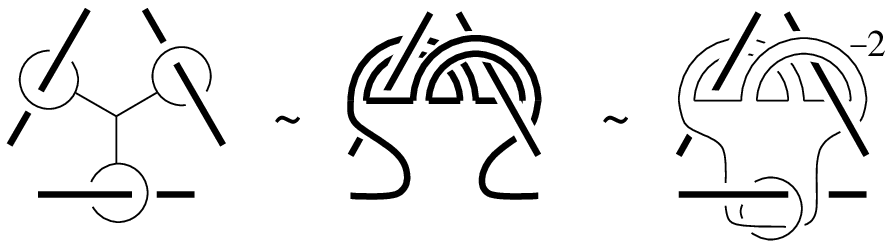,height=1.0in}$$
\end{lem}
% The framings $0$ and $+2$ of the curves in the
% above figure differ from the preffered framings
% by $\pm$.
Therefore the following theorem holds:
\begin{thm}\label{thm_BLeqY}
For each integer $n$ we have
$\cFbl_n(M)=\cFY_n(M)$.
\end{thm}

\subsection{Comparison of $\cFa$ and $\cFY$:
plan of the proof}
The rest of the section is arranged in the following way.

We will first present any Y--link of degree $d$
by a trivial $d$--component Y--link in $S^3$,
together with a trivial unimodular link $O$
which links $T$ in a special way.
We shall then prove an inclusion
$\cFY_{2n}\subset\cFa_{3n}$ by an easy counting
argument.

The opposite inclusion is similar in spirit.
Now we present any $n$--component algebraically
split link by a trivial link $O$, together with a
trivial Y--link $T$ linking $O$ in a special way.
We then consider an appropriately modified
version of some results of Section \ref{sec_graded}.
Building on these results, we prove the inclusion
$\cFa_{3n}\subset\cFY_{2n}\oZ$ using a similar
counting argument.

\subsection{Undoing Y--graphs}\label{sub_untangle}
A link $O$ in $M\sminus G$ {\em laces a clover} $G$,
if $O$ is trivial unimodular, and each of the
(pairwise disjoint) discs bounding its components
intersects $G$ in at most two points, which belong to
the leaves of $G$.
A  Y--link in $M$ is {\em trivial}, if it consists %%%% $G$ deleted %%%%
of Y--graphs standardly embedded in $n$ disjoint balls.
\begin{lem}\label{lem_triv}
Let $T$ be a trivial $n$--component Y--link in $S^3$.
For any $n$--component Y--link $G$ in $S^3$, there
exists a link $O$ in $S^3$ which laces $T$, such that
$[S^3,G]=[S^3_O,T]$.
\end{lem}
\begin{proof}
Any Y--link in $S^3$, in particular $G$, can be made
into a trivial Y--link by framing twists and crossing
changes. Moreover, it suffices to use crossing changes
which involve only the leaves of Y--graphs.
Indeed, instead of a crossing change which involves an
edge of a Y--graph, one can do two subsequent crossing
changes with the neighboring leaf of this graph (by
sliding first the other branch towards the leaf).
Using $K_4$ we can realize each of these framing and
crossing changes by surgery on a trivial unimodular
surgery link, as illustrated in Figure
\ref{fig_kirby_ex}.
%$$\fig{clasp.eps,height=0.60in}$$
The resulting collection of these surgery components
comprises $O$.
\end{proof}
For a Y--link $G$ in a manifold $M$ and a link
$L\subset M\sminus G$ denote by $[M,G,L]$ the
{\em double alternating sum}:
$$[M,L,G]=[[M,L],G]=[[M,G],L]=
\sum_{G'\subset G}\sum_{L'\subset L}
(-1)^{|G'|+|L'|}M_{G'\cup L'}$$
\begin{cor}\label{cor_gener}
The space $\cFY_n$ is generated by all $[S^3,O,T]$,
where $T$ is a trivial Y--link in $S^3$ of degree
at least $n$, and the link $O$ laces $T$.
\end{cor}
\begin{proof}
By Lemma \ref{lem_triv}, $\cFY_n$ is generated by
all $[S^3_O,T]$, with $O$ and $T$ as above.
It remains to notice that
$S^3_O=\sum_{O'\subset O}(-1)^{|O'|}[S^3,O']$
and that any sublink $O'$ of $O$ also laces $T$.
\end{proof}

\begin{thm}\label{YinAS}
For each integer $n$ we have $\cFY_{2n}\subset\cFa_{3n}$.
\end{thm}
\begin{proof}
Let $T$ be a trivial  Y--link in $S^3$ of degree at
least $2n$ and $O$ be an arbitrary link lacing $T$.
In view of Corollary \ref{cor_gener}, it suffices to
prove that $[S^3,O,T]$ belongs to $\cFa_{3n}$.
Suppose that some leaf of $T$ is not linked with $O$,
ie, bounds a disc which does not intersect $O$; then
$[S^3,O,T]=0\in\cFa_{3n}$ by Lemma \ref{lem_disc}.
Otherwise, all (ie, at least $6n$) leaves of $T$
are linked with $O$.
But each component of $O$ is linked with at most two
leaves of $T$; hence the number of components of $O$
is at least $3n$.
Therefore $[S^3,O,T]\in\cFa_{3n}$.
\end{proof}

\subsection{Undoing an AS--link}
Let $L$ be a framed link in $S^3$.
We call $L$ an {\em AS--link}, if it is algebraically
split and unimodular.
In Section \ref{sub_untangle} above we presented any
Y--link by a trivial Y--link $T$, together with a trivial
unimodular link $O$ lacing it.
In this section we shall do the opposite: we present
any AS--link by a trivial unimodular link $O$, together
with a trivial Y--link linking it in a special way,
which we will also call ``lacing".

A clover $G$ in $M\sminus L$ {\em laces a link} $L$,
if each leaf of $G$ either is trivial and links $L$
once (ie, bounds a disc which intersects $L$ in one
point), or is unlinked with $L$ (ie, bounds a surface
which does not intersects $L$).\footnote{In fact, we
will not need this type of leaf, but prefer to
formulate the definition in its full generality}
A pair $(O,G)$ consisting of a trivial link $O$ in $M$
and a Y--link $G$ lacing $O$ is called a {\em lacing pair}.
A lacing pair $(O,G)$ is trivial, if $G$ is trivial.
Surgery on a trivial lacing pair in $S^3$
was called a Borromean surgery in \cite{Ma} and a
$\Delta$--move in \cite{MN}.
It was shown in \cite{Ma, MN}, that one can
pass from a link $L$ in $S^3$ to any other link
with the same linking matrix by surgery on a
trivial Y--link lacing $L$.
Applying this result to AS--links, we deduce:
\begin{lem}\label{lem_trivAS}
Let $O$ be a trivial unimodular $n$--component link
in a $S^3$.
For any AS--link $L$ in $S^3$, there exists a trivial
lacing pair $(O,G)$, such that $[S^3,L]=[S^3_T,O]$.
\end{lem}
\begin{cor}\label{cor_generAS}
The space $\cFa_n$ is generated by all $[S^3,O,T]$,
where $O$ is a trivial unimodular link in $S^3$ with at
least $n$ components and a trivial Y--link $G$ laces $O$.
\end{cor}
\begin{proof}
By Lemma \ref{lem_trivAS}, $\cFa_n$ is generated
by all $[S^3_T,O]$, with $T$ and $O$ as above.
It remains to notice that
$S^3_T=\sum_{T'\subset T}(-1)^{|T'|}[S^3,T']$
and that any Y--sublink $T'$ of $T$ also laces $O$.
\end{proof}
In what follows, we will in fact need only a weaker
version of Corollary \ref{cor_generAS}, in which we
omit the assumption of triviality of the lacing pair.

\subsection{The inclusion $\cFa_{3n}
            \subset\cFY_{2n}\oZ$}\label{sub_ASinY}
We will need a modification of Corollary
\ref{cor_Y3} and Lemma \ref{lem_frame}.
Returning to their proofs, we notice that both
statements can be stated for Y--links lacing a
fixed link in $M$:
\begin{cor}\label{cor_Y3AS}
Fix a link $L$ in a manifold $M$.
Let $G$ be a Y--link of degree $d$ lacing $L$ and
let $G'$ be obtained from $G$ by sliding an edge
of $G$ along a knot $K\subset M\sminus L\sminus G$.
Suppose that for every Y--link $\GG$ of degree
$(d+1)$ lacing $L$ one has $[M,L,\GG]\eq{k}0$.
Then $[M,G']\eq{k}[M,G]$.
\end{cor}
\begin{lem}\label{lem_frameAS}
Fix a link $L$ in a manifold $M$.
Let $G$ be a Y--link of degree $d$ lacing $L$ and
let $G'$ be obtained from $G$ by twisting the
framing of some edge by a half twist.
Suppose that for every Y--link $\GG$ of degree
$(d+1)$ lacing $L$ one has $[M,L,\GG]\eq{k}0$.
Then $[M,G']\eq{k}-[M,G]$
\end{lem}
Moreover, the following version of Lemma
\ref{lem_fragm2} holds:
\begin{lem}\label{lem_specialAS}
Let $(O,G)$ be a lacing pair in a manifold $M$ with
$G$ of degree $d$.
Suppose that some disc $D_i$ bounding a component
of $O$ intersects $G$ in just one point, which
belongs to a leaf of $G$.
Suppose also that for any lacing pair $(O,\GG)$ with
$\GG$ of degree greater than $d$ we have
$[M,O,\GG]\eq{k}0$.
Then $2[M,G]\eq{k}0$.
\end{lem}
\begin{proof}
Proceeding similarly to the proof of Lemma
\ref{lem_fragm2}, we rotate $O$ together with
the trivial leaf linked with it.
This adds a half twist to the framing of the
adjacent edge of $G$, while preserving
the isotopy class of $O$ and $G$.
Thus $[S^3,O,G]\eq{k}-[S^3,O,G]$ by Lemma
\ref{lem_frameAS}.
\end{proof}
We are in a position to prove the second inclusion
theorem.
\begin{thm}\label{thm_ASinY}
For each integer $n$ we have
$\cFa_{3n}\subset\cFY_{2n}\oZ$.
\end{thm}
\begin{proof}
Let $O$ be a lacing pair in $S^3$ with $O$ having at
least $n$ components.
In view of Corollary \ref{cor_generAS}, it suffices
to prove that $[S^3,O,G]\in\cFY_{2n}\oZ$.
We proceed by downward induction on the degree
$d$ of $G$.
If $d\ge 2n$, then obviously $[S^3,O,G]\in\cFY_{2n}$
and the theorem follows.
Suppose now that the inclusion holds for all Y--links
of degree higher than $d$ and let us prove it for a
Y--link $G$ of degree $d$.

By a repeated use of Corollary \ref{cor_Y3AS} we can
reduce the problem to the case when none of the edges
of $G$ pass through the discs $D_i$ bounding the
components of $O$.

If for some $i$ the disc $D_i$ does not intersect
$G$, then $[S^3,O,G]=0\in\cFY_{2n}$ by Lemma
\ref{lem_disc}, and we are done.

If some disc $D_i$ intersects $G$ in exactly one point
belonging to a leaf $l$ of $G$, the statement follows
from Lemma \ref{lem_specialAS} (applicable by the
induction assumption).

We are left with the case when each component of
$O$ is linked with at least two leaves of $G$.
But each leaf of $G$ can be linked with at most
one component of $O$. Therefore, $G$ should have
at least $6n$ leaves, ie, at least $2n$
components.
Hence $[S^3,O,G]\in\cFY_{2n}$, and the theorem follows.
\end{proof}

\begin{rem}
Over the integers, a simplified version of the
above counting argument leads to an inclusion
$[S^3,O,G]\in\cFY_n$ that does not require Lemmas
\ref{lem_frameAS} and \ref{lem_specialAS}.
Indeed, we proceed as in the proof of Theorem
\ref{thm_ASinY} above, first reducing the problem to
the case when no edges of $G$ pass through the discs
$D_i$, and then noticing that if some $D_i$ does not
intersect $G$, then $[S^3,O,G]=0$ by Lemma \ref{lem_disc}.
Otherwise, each of the $3n$ discs $D_i$ intersects
at least one leaf of $G$.
Each leaf of $G$ can intersect at most one disc,
so $G$ should have at least $3n$ leaves, ie, at
least $n$ components, hence $[S^3,O,G]\in\cFY_n$.
\end{rem}

\section{The equivalence of the $\cFb$ and $\cFY$
filtrations}

\subsection{Plan of the proof}
The section is arranged in the following way.

In the first part we shall prove an inclusion
$\cFY_{2n}\subset\cFb_{n}$.
This is done in two steps.
By Corollary \ref{cor_gener} the space $\cFY_d$
is generated by all $[S^3,O,T]$, where $T$ is a
trivial $d$--component Y--link in $S^3$, and the
link $O$ laces $T$.
In the first step, we take $d=2n^2$, construct
a ``good" sublink $B\subset O$ with at least $n$
components, and show that this implies
$\cFY_{2n^2}\subset\cFb_{n}$.
Building on this result, in the second step we
prove a stronger inclusion
$\cFY_{2n}\subset\cFb_{n}$ using the IHX relation.

In the second part we shall prove an opposite
inclusion $\cFb_n\subset\cFY_{2n}\oZ$.
We present any $n$--component boundary link as an
image of a trivial unimodular link $O$ under
surgery on a Y--link $G$, linking $O$ in a special
way, which we will call ``1--lacing''.
Using this presentation together with an appropriately
modified version of Lemma \ref{lem_specialAS}, we prove
that $\cFb_n\subset\cFY_{2n}\oZ$ by an easy counting
argument.

\subsection{Constructing a boundary sublink}
Let $T$ be a trivial $2n^2$ component Y--link in $M$
and let $O$ be a link which laces $T$.
Below we construct a sublink $B\subset O$ with at
least $n$ components, which will be a boundary link
in the manifold, obtained from $M\sminus(O\sminus B)$
by surgery on $T$.
This construction is based on a notion of a good
sublink.

A link $B$ in $M$  is {\em good}, if $B$ laces $T$
and the following conditions hold:
\begin{itemize}
\item No Y--graph of $T$ is linked with $B$ by all
three leaves.
\item If some Y--graph of $T$ is linked with $B$ by
two leaves, then it is linked with just
one component of $B$.
\end{itemize}

\begin{figure}[hbt]
\centerfig{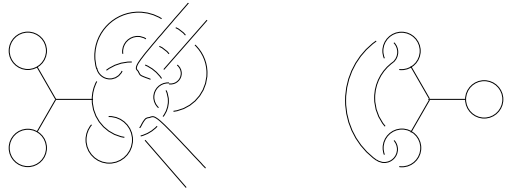,height=1.0in}
\caption{The structure of a good link near a Y--graph}
\label{fig_good}
\end{figure}

The possible structure of $B$ in a neighborhood of a
component of $T$ is shown in Figure \ref{fig_good}.

\begin{lem}\label{lem_bound}
Let $T$ be a trivial Y--link in $M$ and $B$ be a good
sublink.
Then the image $B'$ of $B$ in $M_T$ is a boundary
link in a neighborhood of $T\cup D$, where $D$
is a union of discs bounding the components of $B$.
\end{lem}
\begin{proof}
Cut out a tubular neighborhood $N$ of $T$; this
is a disjoint union of genus 3 handlebodies.
Components of $B$ bound non-intersecting discs in
$M$, which intersect the boundary $\partial N$ by
either one or two circles, which are meridians of $N$.
Perform surgery on $T$ by cutting, twisting
by an appropriate element of the mapping class
group, and gluing back the handlebodies $N$.
This transforms the meridians into some other
curves on $\partial N$.
It suffices to show that for each handlebody
$N_i$ of $N$ these curves on $\partial N_i$ are
bounding inside $N_i$.

Let $m$ be a meridian on $\partial N_i$.
The surgery on this Y--graph corresponds to the
action of an element of the Torelli subgroup
(which acts trivially on $H_1(N_i)$), so the
image of $m$ still bounds a surface $\GS_m$.
Now, recall that $B$ is a good link, so there are
only two possible configurations of the meridians.
Firstly, $\partial N_i$ may contain just two meridians
$m$ and $m'$, corresponding to the same component
$b$ of $B$.
Then, smoothing the intersections of $\GS_{m}$ with
$\GS_{m'}$, we obtain a surface bounding $m\cup m'$.
Secondly, these meridians may appear only on one of
the three handles of $\partial N_i$.
Pick one of the meridians $m$; all other meridians
on $\partial N_i$ may be considered as small
push-offs of $m$ in the normal direction.
Thus their images under the Y--surgery, together
with the corresponding surfaces, can be obtained
by a similar push-offs.
This finishes the proof of the lemma.
\end{proof}

\begin{cor}\label{cor_subset}
Let $T$ be a trivail Y--link in $M$ and $B$ be an
$n$--component good link with the union $D$ of
discs bounding $B$.
Then $[M,L\cup B,T]\in\cFb_n$ for any link $L$ in
$M\sminus(T\cup D)$.
\end{cor}
\begin{proof}
By Lemma \ref{lem_bound}, the image of B under
surgery on $T'\subset T$ bounds in the complement
of $L$. Thus it will bound also in $(M_{T'})_{L'}$,
for any $L'\subset L$.
\end{proof}

\begin{lem}\label{lem_many}
Let $T$ be a trivial $2n^2$--component Y--link in $M$
and let $O$ be a link, which laces $T$ and links all
leaves of $T$.
Then there exists a good sublink $B$ of $O$ with at
least $n$ components.
\end{lem}
We will say that two components of $O$ (respectively,
two leaves of $T$) are {\em neighboring}, if there is
a Y--graph (respectively, a component of $O$) which is
linked with both of them.
\begin{proof}
Suppose that there is a leaf of $T$, which has at
least $n$ neighbors, apart from the other leaves
of the same Y--graph.
Then the corresponding components of $O$ comprise $B$.
It remains to consider the case when each leaf of $T$
has less than $n$ neighbors belonging to other Y--graphs.
Pick an arbitrary component $B_1$ of $O$ and remove
from $O$ all neighboring components of $B_1$.
Repeat this step, each time picking a new component
$B_i$ of the remaining link, until there are no more
components left.
Finally, take $B=\cup_i B_i$.

Let us establish a lower bound for the number of
these steps.
Each $B_i$ is linked with at most two Y--graphs, leaves
of which have altogether at most $6(n-1)-2=6n-8$ other
neighboring leaves.
Thus the removal of $B_i$ and of the link components
neighboring $B_i$ may unlink at most $6n-2$ leaves
of $T$.
In the beginning, $O$ was linked with all $6n^2$
leaves of $G$.
Therefore, the number of steps is at least
$\frac{6n^2}{6n-2}>n$, thus $B$ has more than $n$
components.
\end{proof}

\begin{rem}
In fact, one can assume that the degree of $T$ is
just $6n$, but the construction of $B$ in this case
is significantly more complicated.
\end{rem}

\subsection{The inclusion $\cFY_{2n}\subset\cFb_n$}
Let us start with a weaker inclusion:
\begin{prop}\label{prop_YinB}
For each integer $n$, we have
$\cFY_{2n^2}\subset\cFb_n$.
\end{prop}
\begin{proof}
Let $T$ be a trivial $2n^2$--component Y--link
in $S^3$ and $O$ be an arbitrary link lacing $T$.
In view of Corollary \ref{cor_gener}, it suffices
to prove that $[S^3,O,T]\in\cFb_n$.
If some leaf of $T$ is not linked with $O$, then
$[S^3,O,T]=0\in\cFb_n$.
Otherwise, all leaves of $T$ are linked with $O$
and we can use Lemma \ref{lem_many} to find a good
sublink $B$ of $O$ with at least $n$ components.
It remains to apply Lemma \ref{lem_bound} for $M=S^3$
and $L=O\sminus B$.
\end{proof}

Now we are in a position to prove a stronger result:
\begin{thm}
For each integer $n$, we have $\cFY_{2n}\subset\cFb_n$.
\end{thm}
\begin{proof}
Let $G$ be a clover of degree $d\ge 2n$ in $S^3$.
We proceed by the downward induction on $d$.
If $d\ge 2n^2$, then $[S^3,G]\in\cFb_n$ by
Proposition \ref{prop_YinB}, and the theorem
follows.
Suppose that the statement holds for any clover
of degree higher than $d$, and let us prove it
for a clover $G$ of degree $d$.
Note that by the induction assumption we can use
Theorem \ref{thm_graph}.
Thus it suffices to prove the statement for a
clover $G$ which is a disjoint union of a degree
$2k\le d$ clover $G'$ with no leaves and $d-2k$
copies $G_i$ of a Y--graph with three special leaves.

We call a path in a connected graph {\em maximal},
if it is connected, passes along each edge at most
once, and contains the maximal number of edges.
A path in $G'$ is maximal, if it is maximal in each
of its connected components.
The number $v$ of vertices of $G'$ which do not belong
to a maximal path is called the {\em length-defect}
of $G'$.
If $v$ is positive, pick a vertex of $G'$ which
does not belong to a maximal path, but is connected
to it by an edge.
Applying to this edge the IHX relation of Theorem
\ref{thm_IHX}, we obtain two clovers, each of
which has the length defect $v-1$.
Hence it suffices to prove the theorem for $v=0$.
Notice that if $v=0$, there are at least $k$ edges
of $G'$ without common ends.
Thus, including $G_i$'s, there are at least $d-k$
edges of $G$ without common ends.

Cut all edges of $G'$ by Theorem \ref{thm_edgecut}.
Now, unlink all pairs of newly created leaves and
change the framing of all leaves of $G_i$'s to 0
by $K_4$, as shown below:
$$\fig{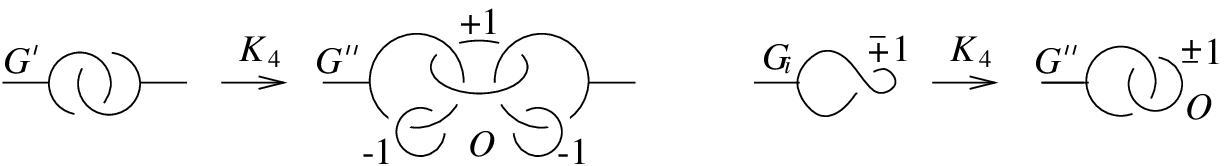,height=0.6in}$$
Denote by $O$ the resulting link and by $T$ the
$d$--component Y--link obtained from $T$.
Clearly the link $O$ laces $T$ and
$$[S^3,G]=\pm[S^3_O,T]=\pm\sum_{O'\subset O}
(-1)^{|O'|}[S^3,O',T].$$
We should show that $[S^3,O',T]\in\cFb_n$ for
each $O'\subset O$.
If some leaf of $T$ is unlinked with $O'$,
then this leaf is trivial and $[S^3,O',T]=0$.
Otherwise, by the construction of $O$, there
are at least $d-k\ge d/2\ge n$ non-neighboring
components of $O'$ (since there were at least
$d-k$ edges of $G$ without common ends).
These components comprise a good sublink of $O'$.
The theorem now follows from Corollary \ref{cor_subset}.
\end{proof}

\subsection{Undoing a boundary link}
Let $B$ be a boundary link in $M$, and fix
a surface $\GS=\cup_i\GS_i$ bounded by $B$.
 We would like to present $B$ by a trivial
unimodular link bounding a collection of discs,
together with a Y--link, lacing it in a rather
special way.
Namely, we will say that a Y--link {\em $b$--laces the
link} $B$, and each of its components intersects $\GS$
in at most one point, which belongs to a trivial leaf.
By the classification of surfaces, we can assume
that a surface $\GS_i$ bounding each component is
an embedding of a connected sum of discs with two
attached bands.
Lemma \ref{lem_blink} shows how such a disc with two
(possibly linked and knotted) bands can be obtained from
a standard disc by surgery on a Y--graph $b$--lacing it:
%\begin{lem}\label{lem_admissible}
%We have:
$$\fig{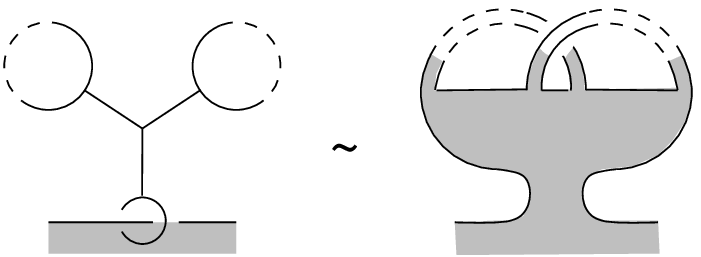,height=0.85in}$$
%\end{lem}
A pair $(O,G)$ consisting of a trivial link $O$ and
a Y--link $G$ in $M$ is called a {\em $b$--lacing pair},
if $G$ $b$--laces $O$ for $\GS_i=D_i$ being discs bounding
the components of $O$.
Exchanging pairs of bands to $b$--lacing Y--graphs
(or vice versa) as above, we obtain
\begin{cor}\label{cor_trivB}
For any $n$--component boundary link $B$ in a manifold
$M$ there exists a $b$--lacing pair $(O,G)$ such that surgery
on $G$ transforms $O$ into $B$.
Conversly, for any $b$--lacing pair $(O,G)$, surgery on
$G$ transforms $O$ into a boundary link.
\end{cor}
We illustrate this construction on an example of
a genus 2 surface bounding a trefoil:
$$\fig{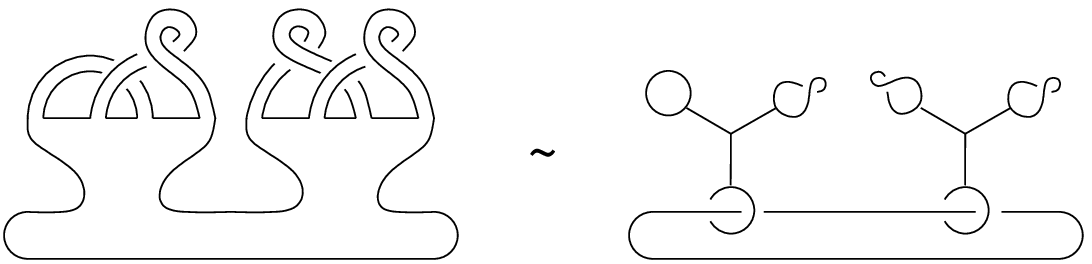,height=0.85in}$$
\begin{rem}\label{rem_plan}
A presentation of a boundary link by a $b$--lacing
pair $(O,G)$ encodes an information about the Seifert
surface $\GS$ of $B$ and the Seifert form.
Indeed, each Y--graph $G_i$ of $G$ corresponds to a
pair of bands, ie, to a handle $h_i$ of $\GS$.
The cores $b_{2i-1}$, $b_{2i}$ of these bands, being
appropriately oriented and closed in the disc in a
standard way, comprise a preferred basis of $H_1(\Sigma)$.
This observation allows one to deduce that the
Seifert matrix is given by $\lk(b_i,b_j)+\Gd_{i+1,j}$,
where $\lk(b_i,b_j)$ is the linking matrix of $2n$
non-trivial leaves of $G$.
We will investigate some applications of this
construction in a future paper.
\end{rem}
\begin{cor}
The space $\cFb_n$ is generated by all $[S^3,O,G]$,
where $(O,G)$ is a $b$--lacing pair in $S^3$ and $O$
has at least $n$ components.
\end{cor}
\begin{proof}
By Corollary \ref{cor_trivB}, $\cFb_n$ is generated
by all $[S^3_G,O]$, with $G$ and $O$ as above.
It remains to notice that
$S^3_G=\sum_{G'\subset G}(-1)^{|G'|}[S^3,G']$,
and that any Y--sublink $G'$ of $G$ also $b$--laces $O$.
\end{proof}

\subsection{The inclusion $\cFb_n\subset\cFY_{2n}$}
Returning to the proof of Lemma \ref{lem_specialAS}
(and Lemmas  \ref{lem_frameAS}, \ref{lem_frame} used
in it), we notice that it may be restated for
$b$--lacing pairs in $M$.
Thus the following modified version of Lemma
\ref{lem_specialAS} holds:
\begin{lem}\label{lem_specialB}
Let $(O,G)$ be a $b$--lacing pair in a manifold $M$
with $G$ of degree $d$.
Suppose that some disc $D_i$ bounding a component
of $O$ intersects $G$ in just one point.
Suppose also that for any $b$--lacing pair $(O,\GG)$
with $\GG$ of degree greater than $d$ we have
$[M,O,\GG]\eq{k}0$.
Then $2[M,G]\eq{k}0$.
\end{lem}
We are ready to prove the last inclusion theorem.
\begin{thm}\label{thm_BinY}
For each integer $n$ we have
$\cFb_n\subset\cFY_{2n}\oZ$.
\end{thm}
\begin{proof}
Let $(O,G)$ be a $b$--lacing pair in $S^3$ such that
$O$ has at least $n$ components.
We will show by downward induction on the degree $d$
of $G$ that $[S^3,O,G]\in\cFY_{2n}\oZ$.
If $d\ge 2n$, then obviously $[S^3,O,G]\in\cFY_{2n}$.
Inductively suppose that the statement holds for any
$b$--lacing pair $(O,\GG)$ in $S^3$ with $\GG$ of degree
greater than $d$.

If for some $i$ the disc $D_i$ bounding a component
of $O$ does not intersect $G$, then obviously
$[S^3,O,G]=0\in\cFY_{2n}$ and we are done.
If for some $i$ the disc $D_i$ intersects $G$ in
exactly one point, the statement follows from Lemma
\ref{lem_specialB} (applicable by the induction
assumption).
We are left with the case when each disc $D_i$
intersects $G$ in at least two points.
Since $G$ $b$--laces $O$, these points should belong
to the leaves of different Y--graphs, which do not
intersect any other disc $D_j$, $j\ne i$.
Therefore, $G$ should have at least $2n$ components,
and the theorem follows.
\end{proof}

\begin{rem}
The inclusion of Theorem \ref{thm_BinY} and its
proof are valid not just for $\Z$HS, but for
arbitrary 3--manifolds.
\end{rem}

\end{document}